\documentclass[11pt,reqno]{amsart}
\usepackage{enumerate}
\usepackage{amsfonts}
\usepackage{amsmath}
\usepackage{amsthm}
\usepackage{latexsym}
\usepackage{multicol}
\usepackage{verbatim}
\usepackage{amscd}
\usepackage{amssymb}
\usepackage{hyperref}
\usepackage{amsmath, amscd,wasysym}
\usepackage{graphicx}
\usepackage{booktabs,caption,fixltx2e}
\usepackage[flushleft]{threeparttable}
\usepackage{tikz}
\usetikzlibrary{chains}
\graphicspath{ {images/} }
\usepackage[all,cmtip,arrow]{xy}
\usepackage{pb-diagram}

\advance\textwidth by 1.2in \advance\oddsidemargin by -.6in \advance\evensidemargin by -.6in

\tikzset{node distance=1.5em, ch/.style={circle,draw,on chain,inner sep=2pt},
fch/.style={circle,fill=black,on chain,inner sep=2pt}, rt/.style={circle,fill=gray,on chain,inner sep=2pt}, fchj/.style={fch,join}, rtj/.style={rt,join},
  chj/.style={ch,join},every path/.style={shorten >=4pt,shorten <=4pt},line width=1pt,baseline=-1ex}

\let\dlabel=\alabel
\newcommand{\dfnode}[2][fchj]{%
\node[#1,label={below:\dlabel{#2}}] {};
}
\newcommand{\drnode}[2][rtj]{%
\node[#1,label={below:\dlabel{#2}}] {};
}

\newcommand{\dnode}[2][chj]{%
\node[#1,label={below:\dlabel{#2}}] {};
}

\newcommand{\dnodenj}[1]{%
\dnode[ch]{#1}
}
\newcommand{\drnodenj}[1]{%
\dnode[rt]{#1}
}

\newcommand{\dfnodenj}[1]{%
\dfnode[fch]{#1}
}

\newcommand{\dnodebr}[1]{%
\node[chj,label={below right:\dlabel{#1}}] {};
}
\newcommand{\drnodebr}[1]{%
\node[rtj,label={below right:\dlabel{#1}}] {};
}

\newcommand{\dfnodebr}[1]{%
\node[fchj,label={below right:\dlabel{#1}}] {};
}

\newcommand{\dydots}{%
\node[chj,draw=none,inner sep=1pt] {\dots};
}

\newtheorem*{cor}{Corollary}
\newtheorem*{lem}{Lemma}
\newtheorem*{prop}{Proposition}

\theoremstyle{definition}
\newtheorem*{defn}{Definition}

\theoremstyle{definition}
\newtheorem*{thm}{Theorem}

\newtheorem*{rem}{Remark}
\newtheorem*{rems}{Remarks}

\newcounter{cnt}
 \makeatletter
\def\mydggeometry{\makeatletter\dg@YGRID=1\dg@XGRID=20\unitlength=0.003pt\makeatother}
\makeatother \theoremstyle{remark}

\numberwithin{equation}{section}

\let\bwdg\bigwedge
\def\bigwedge{{\textstyle\bwdg}}

\newcommand{\thmref}[1]{Theorem~\ref{#1}}

\newcommand{\secref}[1]{Section~\ref{#1}}
\newcommand{\lemref}[1]{Lemma~\ref{#1}}
\newcommand{\propref}[1]{Proposition~\ref{#1}}
\newcommand{\corref}[1]{Corollary~\ref{#1}}
\newcommand{\remref}[1]{Remark~\ref{#1}}

\newcommand{\wt}{\operatorname{wt}}

\newcommand{\nc}{\newcommand}
\newcommand{\rnc}{\renewcommand}

\nc{\cal}{\mathcal} \nc{\goth}{\mathfrak} \rnc{\bold}{\mathbf}
\nc{\fk}{\mathfrak}
\newcommand{\supp}{\operatorname{supp}}

\renewcommand{\Bbb}{\mathbb}
\nc\bomega{{\mbox{\boldmath $\omega$}}}\nc\bOmega{{\mbox{\boldmath $\Omega$}}}
\nc\bchi{{\mbox{\boldmath $\chi$}}}
 \nc\bpsi{{\mbox{\boldmath $\Psi$}}} \nc\bast{{\mbox{\boldmath $\ast$}}}
 \nc\balpha{{\mbox{\boldmath $\alpha$}}}
 \nc\bpi{{\mbox{\boldmath $\pi$}}}
\nc\bsigma{{\mbox{\boldmath $\sigma$}}} \nc\bcN{{\mbox{\boldmath $\cal{N}$}}} \nc\bcm{{\mbox{\boldmath $\cal{M}$}}} \nc\bLambda{{\mbox{\boldmath
$\Lambda$}}} 
\nc\bll{{\mbox{\boldmath$\ell$}}} \nc\bgamma{{\mbox{\boldmath$\gamma$}}}
\renewcommand\i{{\mbox{\boldmath$\iota$}}}
\nc\p{{\mbox{\boldmath$\rho$}}}

\newcommand{\lie}[1]{\mathfrak{#1}}
\makeatletter
\def\section{\def\@secnumfont{\mdseries}\@startsection{section}{1}%
  \z@{.7\linespacing\@plus\linespacing}{.5\linespacing}%
  {\normalfont\scshape\centering}}
\def\subsection{\def\@secnumfont{\bfseries}\@startsection{subsection}{2}%
  {\parindent}{.5\linespacing\@plus.7\linespacing}{-.5em}%
  {\normalfont\bfseries}}
\makeatother

 \nc{\Hom}{\operatorname{Hom}}
  \nc{\td}{\operatorname{\tilde{d}}}\nc{\D}{\operatorname{d}}
  \nc{\mode}{\operatorname{mod}}
\nc{\End}{\operatorname{End}} \nc{\wh}[1]{\widehat{#1}} \nc{\Ext}{\operatorname{Ext}} \nc{\ch}{\text{ch}} \nc{\ev}{\operatorname{ev}}
\nc{\Ob}{\operatorname{Ob}} \nc{\soc}{\operatorname{soc}} \nc{\rad}{\operatorname{rad}} \nc{\head}{\operatorname{head}}
\nc{\A}{\operatorname{\bold{b}}}
\def\Im{\operatorname{Im}}
\def\ker{\operatorname{Ker}}

 \nc{\Cal}{\cal} \nc{\Xp}[1]{X^+(#1)} \nc{\Xm}[1]{X^-(#1)}
\nc{\on}{\operatorname} \nc{\Z}{{\mathbb{Z}}} \nc{\J}{{\cal J}} \nc{\C}{{\mathbb{C}}} \nc{\Q}{{\bold Q}} \nc{\R}{{\mathbb{R}}}
 
\nc{\N}{{\Bbb N}} \nc\boa{\bold a} \nc\bob{\bold b} \nc\boc{\bold c} \nc\bod{\bold d} \nc\boe{\bold e} \nc\bof{\bold f} \nc\bog{\bold g}
\nc\boh{\bold h} \nc\boi{\bold i} \nc\boj{\bold j} \nc\bok{\bold k} \nc\bol{\bold l} \nc\bom{\bold m} \nc\bon{\bold n} \nc\boo{\bold o}
\nc\bop{\bold p} \nc\boq{\bold q} \nc\bor{\bold r} \nc\bos{\bold s} \nc\boT{\bold t} \nc\boF{\bold F} \nc\bou{\bold u} \nc\bov{\bold v}
\nc\bow{\bold w} \nc\boz{\bold z} \nc\boy{\bold y} \nc\ba{\bold A} \nc\bb{\bold B} \nc\bc{\bold C} \nc\bd{\bold D} \nc\be{\bold E} \nc\bg{\bold
G} \nc\bh{\bold H} \nc\bi{\bold I} \nc\bj{\bold J} \nc\bk{\bold K} \nc\bl{\bold L} \nc\bm{\bold M} \nc\bn{\bold N} \nc\bo{\bold O} \nc\bp{\bold
P} \nc\bq{\bold Q} \nc\br{\bold R} \nc\bs{\bold S} \nc\bt{\bold T} \nc\bu{\bold U} \nc\bv{\bold V} \nc\bw{\bold W} \nc\bz{\bold Z} \nc\bx{\bold
x} \nc\KR{\bold{KR}} \nc\rk{\bold{rk}} \nc\het{\text{ht }}

\nc\fn{{fin}}  \nc\af{{aff}}  \nc\tr{{tor}} \nc\btilde{\bold{\tilde{\bold{H}}}}

\nc{\mpp}{\rotatebox[origin=c]{180}{\pm}}

\nc\eps{\epsilon}
\nc\toa{\tilde a} \nc\tob{\tilde b} \nc\toc{\tilde c} \nc\tod{\tilde d} \nc\toe{\tilde e} \nc\tof{\tilde f} \nc\tog{\tilde g} \nc\toh{\tilde h}
\nc\toi{\tilde i} \nc\toj{\tilde j} \nc\tok{\tilde k} \nc\tol{\tilde l} \nc\tom{\tilde m}  \nc\ton{\tilde n} \nc\too{\tilde o} \nc\toq{\tilde q}
\nc\tor{\tilde r} \nc\tos{\tilde s} \nc\toT{\tilde t} \nc\tou{\tilde u} \nc\tov{\tilde v} \nc\tow{\tilde w} \nc\toz{\tilde z}
\begin{document}

\title[ Integrable Representations of the Toroidal Lie algebras]{Spectral Characters of a class of  Integrable Representations of Toroidal Lie Algebras}
\author{Tanusree Khandai}

\address{ Indian Institute of Science Education and Research, Mohali, Punjab, India
\newline *School of Mathematical Sciences\\ NISER \\ HBNI, Jatni\\ Odisha-752050\\India}
\email{tanusree@iisermohali.ac.in}
\begin{abstract}
  In this paper we study the subcategory of finite-length objects of the category of positive level integrable representations of a toroidal Lie algebra.
   The main goal is to characterize the blocks of the category. In the cases
 when the underlying finite type Lie algebra  associated with the toroidal Lie algebra is simply-laced, we are able to give a parametrization for the 
blocks. 

\end{abstract}

\maketitle

\section{Introduction}

Toroidal Lie algebras are generalizations of affine Kac-Moody Lie algebras. 
In \cite{MRY,MR} a $k$-toroidal Lie algebra $\cal T_k(\lie g)$ was defined as the  universal central extension of the Lie algebra of polynomial maps from
$(\mathbb{C}^\ast)^k$ to a finite-type Kac-Moody Lie algebra $\lie g_\fn$, for $k\in \mathbb N$.  Ever since the structure and representation theory of these Lie algebras have been extensively studied. In contrast to the one-dimensional center of an affine Kac-Moody Lie algebra  $\lie g_{aff}$, the $k$-toroidal Lie algebras have a $\Z^k$-graded infinite-dimensional center and this makes the study of their representations more interesting. 
The representation theory of $\lie g_{aff}$ has been a subject of interest for
the past three decades\cite{C1,CG1,CM,CPloop,CP,Kac,R2,VV}. In recent years, many of the methods that were developed to study the latter have 
been extended to study integrable representations of toroidal Lie algebras and its quotients by
central ideals of finite co-dimension \cite{CFK,CL,FL,Kh,KL,NS,S,R3,R05,RSF}. In this paper we continue with the project.

Let $\cal I_{fin}$ be the category of integrable $\cal T_k(\lie g)$-modules with 
finite-dimensional weight spaces and let $\cal I_{fin}^\bast$ be the full subcategory  of $\cal I_\fn$ consisting of 
$\cal T_k(\lie g)$-modules on which the central elements act non-trivially. The simple objects of the category $\cal I_\fn$ have been classified in 
\cite{R3, Kh}. It is however known that the category $\cal I_{fin}^\bast$ is not semisimple. One is therefore interested in the descriptions of the blocks of the category. 
In this paper we study the subcategory of finite-length objects in $\cal I_{fin}^\bast$ and in certain cases, give a parametrization for the blocks of this subcategory.


The structure of the category $\cal F$ of finite-dimensional representations  of $\lie g_{aff}$ was studied in \cite{CM}. Defining an equivalence relation on the objects of $\cal F$, it was proved that the category $\cal F$ can be decomposed into blocks that are parameterized by finitely supported functions from $\mathbb C^\ast$ to $\Gamma$, the quotient of the weight lattice of
$\lie g_{fin}$ by the root lattice of $\lie {g}_{fin}$.  Using the exactness of the loop functor $\cal L: \cal F \rightarrow \cal I_\fn$
(which maps a finite-dimensional $\lie g_{aff}$-module $V$ to the integrable
$\lie g_{aff}$-module $V\otimes \C[t^{\pm1}]$ in $\cal I_{\fn}$) it was proved in 
\cite{CG1}, that the category of graded level zero integrable representations
with finite-dimensional weight spaces of $\lie g_{aff}$  can similarly be
decomposed into blocks and these blocks are parametrized by orbits for a
natural action of the group $\C^\ast$ on the set of finitely supported
functions from $\C^\ast$  to $\Gamma$.

The first extension groups for finite-dimensional irreducible representations
of generalized current Lie algebras, twisted current algebra and equivariant
map algebras which include the multiloop Lie algebras was studied in
\cite{Ko, AL} and \cite{NS} respectively.  While, by using the results of 
\cite{NS}, the techniques of \cite{CM,CG1}
can be extended verbatim to obtain block decomposition of the subcategory 
$\cal I_{fin}^{(\bold{0})}$ of $\cal I_{fin}$ on which the center acts trivially, 
they do not help in determining the structure of $\cal I_{fin}^\bast$. One of the
 main problems that arise is the fact that the indecomposable modules in 
$\cal I_{fin}^\bast$ do not have finite-length (or pseudo finite-length) 
property. As a first step towards the study of the structure of the category 
$\cal I_{fin}^\bast$,  we therefore restrict our attention to the
subcategory $\cal J_{int}^+$ of  $\cal I_\fn^\bast$ of finite-length objects of 
positive level. 


Let  $\Pi$ be the set of finitely supported functions from $(\mathbb C^\ast)^{k-1} $ to $P_{aff}^+$, the set of dominant integral weights of $\lie g_{aff}$. It is known from \cite{R3,Kh}  that the simple objects in 
$\cal I_{fin}^\bast$ are parametrized by the set of tuples $\{(\pi, \bos): \pi\in \Pi, \bos\in \Z_0^k\}.$
 Let $X_\pi^\bos$ be the irreducible $\cal T_k(\lie g)$-module corresponding to a pair $(\pi,\bos)\in \Pi\times \Z^k_0$. Associating with every $\pi\in \Pi$ a function $\xi_\pi: (\C^\ast)^{k-1} \rightarrow \Z\times \Gamma$, and considering the natural action of $(C^\ast)^{k-1}$ on the set $\Xi= \{\xi_\pi: \pi\in \Pi\}$, we show that, if   $X_\pi^\bos$ and $X_{\pi'}^{\bos'}$ are irreducible sub-quotients of an indecomposable $\cal T_k(\lie g)$-module of finite-length, then $\xi_{\pi}= \bob.\xi_{\pi'}$ for some $\bob\in (\C^\ast)^{k-1}$. The converse however does not hold in general. 

Using results (from \cite{A1,A2}) on 
the irreducibility of the tensor product of a highest weight integrable 
$\lie g_{aff}$-module and the loop modules for $\lie g_{aff}$,
we classify the functions in $\Pi$ into type $\bi$ and type $\bi\bi$. We prove that if $\lie g_{fin}$ is simply-laced and $\pi,\pi'\in \Pi$ are two functions of type $\bi\bi$  with $\xi_{\pi}=\bob.\xi_{\pi'}$ for some $\bob \in (\C^\ast)^{k-1}$, then  there exists a sequence $X_{\pi_1}^{\bor_1}=U_0,U_1,\cdots,U_r=X_{\pi_2}^{\bor_2}$ of indecomposable $\cal T_k(\lie g)$-modules of finite length in $\cal I_{fin}^\bast$ such that $\Hom_{\cal T_k(\lie g)}(U_i,U_{i+1})\ne 0$ or $\Hom_{\cal T_k(\lie g)}(U_{i+1},U_i)\ne 0$ for $0\leq i\leq r-1$.
Further, if $V$ is an indecomposable $\cal T_k(\lie g)$-module in $\cal I_{fin}^\bast$ having an irreducible sub-quotient isomorphic to $X_\pi^\bos$ for a function $(\pi,\bos)\in \Pi\times \Z_0^{k-1}$ with $\pi$ of type $\bi$, then every irreducible constituent of $V$ is isomorphic to $X_\pi^\bos$. These together lead us towards a parametrization of the blocks in the category $\cal J_{int}^+$ in the case when $\lie g_{fin}$ is of type $A_n,D_n,E_6, E_7$ or $E_8$. 

The paper is organized as follows. After setting the notations for
the paper in \secref{prelim}, the basic properties of the  integrable
representations of toroidal Lie  algebras are recalled in \secref{I.fin.ast}.
In \secref{J.fin}, a restricted subcategory $\cal J_{int}^\pm$ of the category
of $\cal I_\fn^\bast$ is introduced and properties of the simple objects in $\cal J_{int}^+$ are listed. In \secref{Weyl}, the properties of  Weyl modules of $\cal T_k(\lie g)$ proved and results from \cite{RSF} are recalled. Finally, in \secref{BLK} and  \secref{Proof.Prop}, the main results of the paper are stated and proved. Our result  gives a complete parametrization of the blocks in $\cal J_{int}^+$ in the cases when $\lie g_{fin}$ is of type $A_n$, $D_n$, $E_6, E_7$ and $E_8$.

\noindent{\bf{Acknowledgement}}{: I sincerely thank the reviewer for his comments and 
suggestions, which significantly improved the paper.}

\section[Preliminaries]{Preliminaries}\label{prelim}
 \noindent  In this section we fix the notations for the paper and recall the explicit realization of $k$-toroidal Lie algebras from \cite{R2,MR}.
 
\subsection{}\label{section1.1}Throughout the paper $\C,\R$, $\Z$ and $\N$ shall denote the field of complex numbers, real numbers, the set of integers and the set of natural numbers. For a commutative associative algebra $\bold A$,  the set of maximal ideals of $\bold{A}$ shall be denoted by $\max \bold{A}$ and for a Lie algebra $\lie a$ the universal enveloping algebra of $\lie a$ shall be denoted by $\cal U(\lie a)$. For $k\in \N$, a $k$-tuple of integers $(m_1,\cdots,m_k)$ shall be denoted by $\bom$ and let $\Z^k_0= \{\bom = (m_1,\cdots,m_k)\in \Z^k\ |\ m_1=0 \}$. 

\subsection{}  Let $\fk{g}_{fin}$ be a finite-dimensional simple Lie algebra of rank $n$, $\fk{h}_{fin}$ a Cartan subalgebra of $\fk{g}_{fin}$ and $R_\fn$ the root system of $\fk{g}_\fn$ with respect to $\fk{h}_\fn$. Let  $\{\alpha_i: 1\leq i\leq n\}$ (respectively $\{\omega_i: 1\leq i\leq n\}$) be a set of simple roots (respectively fundamental weights ) of $\lie g_\fn$ with respect to $\lie h_\fn$, $R_\fn^+$ (respectively $Q_{fin}$, $P_{fin}$) be the corresponding set of positive roots (respectively root lattice and weight lattice) and  let $\theta$(respectively $\theta_s$) be the highest root ( highest short root) of $R_\fn^+$ if $\lie g$ is simply-laced (respectively not simply-laced).  Let $( .\,| .\,)$ be the non–degenerate symmetric bilinear form on $\lie h_{fin}^\ast$ normalized so that the square length of a long root is 2.
 Let $Q_{fin}^+$ and $P_{fin}^+$ be the $\Z_+$ span of the simple roots and fundamental weights of $(\lie g_\fn,\lie h_\fn)$  and let $W_{fin}$ be the Weyl group of $\lie g_{fin}$.

For $\alpha \in R_\fn^\pm$, let $\lie g_\fn^{\pm\alpha}$ denote the corresponding root space. Let $x_\alpha^\pm \in \lie g_\fn^{\pm\alpha}$ and $\alpha^\vee\in \lie h_\fn$ be  fixed elements such that $\alpha^\vee = [x_\alpha^+,x_\alpha^-]$ and $[\alpha^\vee, x_\alpha^\pm] = \pm 2 x_\alpha^\pm$.
 
 Let $\Gamma = P_{fin}/Q_{fin}$. It is well known that $\Gamma$ is a finite abelian group and the non-zero elements in $\Gamma$ are of the form $\omega_i\mod Q_\fn$ where $\omega_i$ is a fundamental weight of $\lie g_\fn$ such that $\omega_i(\alpha) =0$ or $1$ for all $\alpha\in R_\fn^+.$ 
Let \begin{equation} \label{J.0} 
J_0 = \left\lbrace\begin{array}{lr}
\{1,2,\cdots,n\} & \text{if}\ \lie g \ \text{is of type}\ A_n\\ 
\{n\} & \text{if}\ \lie g \ \text{is of type}\ B_n\\
\{1\} & \text{if}\ \lie g \ \text{is of type}\ C_n\\
\{1,n-1,n\} & \text{if}\ \lie g \ \text{is of type}\ D_n\\
\{1,6\} & \text{if}\ \lie g \ \text{is of type}\ E_6\\
\{7\} & \text{if}\ \lie g \ \text{is of type}\ E_7
\end{array}\right.\end{equation}
Using the labelling of the nodes as in \cite[Plate I-IX]{Bou}, by \cite[Section 13, Exercise 13]{Humphreys}, we have 
$$\Gamma = \{\omega_i \mod Q_\fn : i\in J_0\}\cup\{0\},$$ when $\lie g_\fn$ of type $A_n, B_n, C_n, D_n, E_6, E_7$,  and $\Gamma = \{0\}$, when $\lie g_\fn$ of type $E_8, F_4, G_2$. 
\label{J0}

\vspace{.25cm}
For $\lambda\in P_\fn^+$, let $V(\lambda)$ denote the cyclic $\lie g_\fn$-module generated by a weight vector $v_\lambda$ with defining relations:
$$ x_\alpha^+.v_\lambda =0,\ \forall\ \alpha\in R_\fn^+, \hspace{.5cm} h.v_\lambda =\lambda(h)v_\lambda, \ \forall\ h\in\lie h_\fn,\hspace{.5cm} (x_\alpha^-)^{\lambda(\alpha)+1}.v_\lambda=0, \ \forall\ \alpha\in R^+_\fn.$$ 
It is well known that $V(\lambda)$ is an irreducible finite-dimensional $\lie g_\fn$-module with highest weight $\lambda$ and any irreducible finite-dimensional $\lie g_\fn$-module is isomorphic to $V(\mu)$ for $\mu\in P^+_\fn.$

\vspace{.35cm}

\subsection{} For a positive integer $k$, let $\C[t_1^{\pm 1},\cdots,t_k^{\pm1}]
$ be the Laurent polynomial ring in $k$ commuting variables $t_1,\cdots,t_k$ and for $\bom=(m_1,\cdots,m_k)\in \Z^k$ , let $t^\bom$ denote the element
$t_1^{m_1}\cdots t_k^{m_k}$ in $\C[t_1^{\pm1},\cdots,t_k^{\pm1}]$. Let
$L_k(\fk{g}) = \fk{g}_\fn\otimes \C[t_1^{\pm 1},t_2^{\pm 1},\cdots, t_k^{\pm 1}]$ and $\cal{Z}_k=\bOmega_k/dL_k$ be the space of K$\ddot{\text{a}}$hler differentials  spanned by the set of vectors $\{ t^\bom K_i,\  ~ \bom\in \Z^k, ~ 1\leq i\leq k\}$ together with the  relation, 
$\sum_{i=1}^k r_i t^{\bor} K_i = 0$.
Let $d_i:({L}_k(\fk{g})\oplus \cal Z_k)\rightarrow ({L}_k(\fk{g})\oplus \cal Z_k)$, $1\leq i\leq k$ be the $k$ derivations on $L_k(\fk{g})\oplus \cal Z_k$ given by: 
\begin{eqnarray}\label{derivation1}
d_i(x\otimes t^\bom)= m_i x\otimes t^\bom, \hspace{.35cm} 
d_i(t^\bom K_j)=m_i t^\bom K_j & \hspace{.15cm} \forall\ \, x\in \lie g_{fin},\, \, \bom\in \Z^k,\, 1\leq i,j\leq k,\,
\end{eqnarray} and let $D_k$ be the $\C$ linear span of the derivations $d_1,d_2\cdots,d_k$. The $k$-toroidal Lie algebra associated to a simple Lie algebra $\fk{g}_\fn$ is the vector space $\cal T_k(\lie g) = L_k(\lie g)\oplus \cal Z_k\oplus D_k$ on which the Lie bracket is defined by \eqref{derivation1} and the following
relations:
\begin{equation}\label{central}
[x\otimes t^\bom, y\otimes t^\bos]=[x,y]\otimes t^{\bom+\bos} +\overline{t^\bos(Dt^{\bom})} (x|y),\hspace{.35cm}[x\otimes t^\bom, \omega ] = 0\hspace{.35cm} [\omega, \omega']=0,
\end{equation} where $x,y\in \fk{g}$, $\bom,\bos\in \Z^k$, $\omega, \omega'\in \cal{Z}_k$ and $\overline{t^\bos(Dt^\bom)} = \sum_{i=1}^k m_it^{\bom+\bos}K_i$. Let $\cal Z_0$ is the subspace of $\cal Z_k$ spanned by the central elements $\{K_i: 1\leq i\leq k\}$ and let $\lie h_{tor}:= \lie h_\fn \oplus \cal Z_0\oplus D_k.$ In order to identify $\lie h_\fn^*$ with a subspace of $\lie h_{tor}^*$, an element $\lambda\in \fk{h}_\fn^*$ is extended to an element of $\lie h_{tor}^*$ by setting $\lambda(c)=0 =\lambda(d_i)=0,$ for all $c\in \cal{Z}_0,  1\leq i\leq k.$
For $1\leq i\leq k$, define $\delta_i \in \lie h_{tor}^*$ by
$ \delta_i|_{\fk{h}_\fn+\cal{Z}_0}=0,$ $\delta_i(d_j) =\delta_{ij},$ for $1\leq j \leq k.$ Given $\bom=(m_1,\cdots,m_n)\in \Z^k$, set $\delta_\bom = \sum_{i=1}^k m_i\delta_i$ and let  $$R_\tr^{re}= \{\alpha+\delta_\bom : \alpha\in R_\fn, \bom\in \Z^k \}, \hspace{.5cm}  \hspace{.25cm} R_\tr^{im} = \{\delta_\bom : \bom\in\Z^k-\{\bold{0}\} \}. $$ $R_\tr^{re}$ and $R_\tr^{im}$ are respectively the set of real and imaginary roots of $\cal T_k(\lie g)$ and $R_\tr := R_\tr^{re}\cup R_\tr^{im}$ is the set of all roots of $\cal T_k(\lie g)$ with respect to $\lie h_{tor}.$ The root vector corresponding to a real root $\alpha+\delta_\bom$ is of the form $x_\alpha\otimes t^\bom$ and the root vectors corresponding to an imaginary root $\delta_\bom$ are of the form $h\otimes t^\bom$ with $h\in \lie h_\fn$. Set $\alpha_{n+i}:=\delta_i-\theta$, for $ i=1,\cdots,k$.

Given $\gamma=\alpha+\delta_\bom\in R_\tr^{re}$, with $\alpha\in R_\fn^+$ and $\bom\in \Z^k$, let $\gamma^\vee=\alpha^\vee+\frac{2}{(\alpha|\alpha)}\sum_i m_i K_i.$ With the given Lie bracket operation on $\cal T_k(\lie g)$, it is easy to check that the subalgebra $\lie {sl}_2(\gamma)$ of $\cal T_k(\lie g)$ spanned by the set of elements $\{x_\alpha^+\otimes t^\bom, x_\alpha^-\otimes t^{-\bom}, \gamma^\vee\}$ is isomorphic to $\lie {sl}_2(\C)$.

\label{roots}

\subsection{} \label{affine}
\vspace{.1cm}

\begin{table}
\begin{align*}
&A_1^{(1)} &&
\begin{tikzpicture}[start chain]
\drnodenj{1}
\dfnodenj{1}
\path (chain-1) -- node[anchor=mid] {\(\Longleftrightarrow\)} (chain-2);
\end{tikzpicture} && \\
&A_l^{(1)} (l \ge 2) &&
\begin{tikzpicture}[start chain,node distance=1ex and 2em]
\dfnode{1}
\dfnode{1}
\dydots
\dfnode{1}
\dfnode{1}
\begin{scope}[start chain=br going above]
\chainin(chain-3);
\node[rt,join=with chain-1,join=with chain-5,label={[inner sep=1pt]10:\(1\)}] {};
\end{scope}
\end{tikzpicture}
\\
&B_l^{(1)} (l \ge 3) &&
\begin{tikzpicture}
\begin{scope}[start chain]
\drnode{1}
\dnode{2}
\dnode{2}
\dydots
\dnode{2}
\dfnodenj{1}
\end{scope}
\begin{scope}[start chain=br going above]
\chainin(chain-2);
\dnodebr{1}
\end{scope}
\path (chain-5) -- node{\(\Rightarrow\)} (chain-6);
\end{tikzpicture}
\\
&C_l^{(1)} (l \ge 2) &&
\begin{tikzpicture}[start chain]
\drnodenj{1}
\dfnodenj{1}
\dydots
\dnode{1}
\dnodenj{1}
\path (chain-1) -- node{\(\Rightarrow\)} (chain-2);
\path (chain-4) -- node{\(\Leftarrow\)} (chain-5);
\end{tikzpicture}
\\
&D_l^{(1)} (l \ge 4) &&
\begin{tikzpicture}
\begin{scope}[start chain]
\drnode{1}
\dnode{2}
\dnode{2}
\dydots
\dnode{2}
\dfnode{1}
\end{scope}
\begin{scope}[start chain=br going above]
\chainin(chain-2);
\dfnodebr{1};
\end{scope}
\begin{scope}[start chain=br2 going above]
\chainin(chain-5);
\dfnodebr{1};
\end{scope}
\end{tikzpicture}
\\
&G_2^{(1)} &&
\begin{tikzpicture}[start chain]
\drnode{1}
\dnode{2}
\dnodenj{1}
\path (chain-2) -- node{\(\Rrightarrow\)} (chain-3);
\end{tikzpicture}
\\
&F_4^{(1)} &&
\begin{tikzpicture}[start chain]
\drnode{1}
\dnode{2}
\dnode{3}
\dnodenj{2}
\dnode{1}
\path (chain-3) -- node[anchor=mid]{\(\Rightarrow\)} (chain-4);
\end{tikzpicture}
\\
&E_6^{(1)} &&
\begin{tikzpicture}
  \begin{scope}[start chain]
    \dfnode{1}
\foreach \dyi in {2,3,2} {
\dnode{\dyi}
}
\dfnode{1}
\end{scope}
\begin{scope}[start chain=br going above]
\chainin(chain-3);
\dnodebr{2}
\drnodebr{1}
\end{scope}
\end{tikzpicture}
\\
&E_7^{(1)} &&
\begin{tikzpicture}
  \begin{scope}[start chain]
    \drnode{1}
\foreach \dyi in {2,3,4,3,2} {
\dnode{\dyi}
}
\dfnode{1}
\end{scope}
\begin{scope}[start chain=br going above]
\chainin(chain-4);
\dnodebr{2}
\end{scope}
\end{tikzpicture}
\\
&E_8^{(1)} &&
\begin{tikzpicture}
\begin{scope}[start chain]
\drnode{1}
  \foreach \dyi in {2,3,4,5,6,4,2} {
\dnode{\dyi}
}
\end{scope}
\begin{scope}[start chain=br going above]
\chainin(chain-6);
\dnodebr{3}
\end{scope}
\end{tikzpicture}
\end{align*}

\caption{ Dynkin Diagrams of Non-twisted Affine Kac-Moody Lie algebras}\begin{tablenotes}\small\label{table:affine} \item{In these diagrams the gray nodes correspond to the root $\alpha_{n+1}$ and the remaining nodes are enumerated as in \cite[Plate I-IX]{Bou}. The numerical labels given here correspond to the number $a_i^\vee$( Refer to \cite[Section 6.1, Table Aff 1, Table Aff 2, Table Aff 3]{Kac}) and the blackened nodes correspond to those contained in $J_0$.}
\end{tablenotes}
\end{table}

For $k=1$, the Lie algebra $\cal T_1(\lie g)$ is an affine Kac-Moody Lie algebra and we denote it by $\lie g_{aff}$. Explicitly $\lie g_{aff}=\lie g_\fn\otimes \C[t_1^{\pm1}]\oplus \C K_1 \oplus \C d_1.$
Owing to the natural ordering in $\Z $, the set of real and imaginary roots of
$\lie g_{aff}$ can be partitioned as follows:
$$\begin{array}{cc}
R^{ {re}^{\pm}}_{aff}= \{\alpha+m\delta_1 : \alpha\in R_\fn, m\in\Z_\pm \backslash\{0\} \}\cup R_\fn^\pm, \hspace{.35cm} & \hspace{.35cm}
R^{{im}^\pm}_{aff} = \{m\delta_1 : m\in\Z_\pm \backslash\{0\} \}.\end{array}
$$The set $R_{aff}^+ = R_{aff}^{{re}^+} \cup R_{aff}^{{ im}^+}$ (respectively  $R_{aff}^- = R_{aff}^{{re}^-} \cup R_{aff}^{{ im}^-}$ ) is called the set of positive (respectively negative) roots of $\lie g_{aff}$ and $R_{aff} = R_{aff}^+ \cup R_{aff}^-$ is the set of roots of $\lie g_{aff}$. We denote by $x_{\alpha}^+\otimes t_1^m$, $x_\alpha^-\otimes t_1^m$, $h\otimes t^m, h\in \lie h_{fin}$ the root vectors corresponding to the roots $\alpha+m\delta_1$, $-\alpha+m\delta_1$ and $m\delta_1$ for $\alpha\in R_{fin}^+$ and $m\in \Z.$ Denoting the root space of $\lie g_{aff}$ corresponding to a root $\gamma\in R_{aff}$ by $\lie g_{aff}^{\gamma}$, set
$\lie n_{aff}^\pm = \underset{\gamma\in R_{aff}^\pm}{\bigoplus}(\lie g_{aff}^\gamma)$  and $\lie h_{aff}=\lie h_\fn\oplus \C K_1 \oplus \C d_1.$ Then one has the decomposition
$\lie g_{aff} = \lie n^-_{aff}\oplus \lie h_{aff}\oplus \lie n^+_{aff}.$

 The
set of simple roots $\Delta_{aff}$ and  coroots $\Delta_{aff}^\vee$ of
$\lie g_{aff}$ are respectively given by
$$\begin{array}{lll}
\Delta_{aff} &= &\{\alpha_1,\cdots,\alpha_n,\alpha_{n+1}=\delta_1 -\theta\},\\ 
 \Delta_{aff}^\vee &= &\{\alpha_1^\vee,\cdots,\alpha_n^\vee,\alpha_{n+1}^\vee=K_1-\theta^\vee\}.\end{array}$$
Let $Q_{aff}$(respectively $Q_{aff}^\vee$) be the root lattice (respectively
coroot lattice) for $\lie g_{aff}$. Let $\Lambda_i$ ($i=1,\cdots,n,n+1$) be the
fundamental weights of $\lie g_{aff}$,  where  \begin{equation} \label{Lambda.i} \Lambda_i = a_i^\vee\Lambda_{n+1}+\omega_i, \hspace{.35cm} \text{for each}\ 1\leq i\leq n,\end{equation} and $\Lambda_i(\alpha_j^\vee) =\delta_{ij},$ for $1\leq j\leq n+1$ and
$\Lambda_i(d_1)=0.$ Here $a_i^\vee$ is the integer labelling the $i^{th}$-node in the Dynkin diagrams given in Table \ref{table:affine}. Thus, 
$\hspace{.5cm}\lie h_{aff}^* = \lie h_\fn^* \oplus \C\delta_1 \oplus \C\Lambda_{n+1},$  and  an element
$\lambda$ in $\lie h_{aff}^*$ can be uniquely written as
$$\lambda = \lambda(K_1)\Lambda_{n+1} +\lambda\vert_{\lie h_\fn} +
(\lambda|\Lambda_{n+1})\delta_1,$$ where the inner product $(.|.)$ on $\lie h_{aff}^\ast$  is as defined in \cite[$\S$ 6.2]{Kac} and extends the inner product on $\lie h_{fin}$.
Let $P_{aff} = \sum_{i=1}^{n+1} \Z\Lambda_i +\C \delta_1,$ (respectively $P_{aff}^+ = \sum_{i=1}^{n+1} \Z_+\Lambda_i +\C \delta_1$) be the set of integral
weights(respectively dominant integral weights) of $\lie g_{aff}$.  Let $\succeq$ be the partial order on $P_{aff}$ defined 
by $\lambda\succeq \mu$ if $\lambda,\mu\in P_{aff}$ are such that $\lambda-\mu\in \sum_{i=1}^{n+1} \Z_+ \alpha_{i}$. Given $\lambda,\mu\in P_{aff}$, we shall write $\lambda\succ\mu$ whenever $\lambda\succeq\mu$ but $\lambda\neq \mu$.

Let  $\omega: \lie g_{aff} \rightarrow \lie g_{aff}$ be the Cartan involution on
$\lie g_{aff}$ given by :
\begin{equation} \begin{array}{lll}
x_\alpha^+\otimes t_1^{m} & \mapsto &-x_\alpha^-\otimes t_1^{-m},\\
x_\alpha^-\otimes t_1^{m} & \mapsto  &-x_\alpha^+\otimes t_1^{-m},\\
h \otimes t_1^{m} & \mapsto & -h\otimes t^{-m},\end{array}
\end{equation} for $\alpha\in R_\fn^+$, $m\in \Z$ and $h\in \lie h_\fn$.

\vspace{.15cm}

An $\lie h_\af$-diagonalizable module $V$ of $\lie g_\af$ is said to be integrable if the root vectors corresponding to the real roots of $\lie g_\af$ are locally nilpotent on $V$. An integral $\lie g_\af$-module is said to be of positive (respectively negative) level if the central element $K_1$ acts on $V$ by a positive (respectively negative) integer. Given $\lambda\in P^+_{aff}$, let $X(\lambda)$ be the irreducible integrable $\lie g_{aff}$-module with highest weight $\lambda$ and highest weight vector 
$v_\lambda$. As $X(\lambda)$ is integrable we can write it as a direct sum of 
its weight spaces :
$$X(\lambda) = \bigoplus_{\nu\in P_{aff}} X(\lambda)_{\nu}, \hspace{.5cm} \text{where}\  X(\lambda)_\nu=\{ v\in X(\lambda): h.v=\nu(h)v,\ \forall \ h\in \lie h_{aff}\}.$$ Let $P_{aff}(\lambda) =\{\nu\in P_{aff}: X(\lambda)_\nu \neq 0\}$ and $X^\ast(\lambda) := \bigoplus_{\nu\in P_{aff}} (X(\lambda)_\nu)^\ast\subset X(\lambda)^\ast.$ It was shown in \cite{Kac} that the restricted dual $X^\ast(\lambda)$ is an irreducible integrable $\lie g_{aff}$-submodule of $X(\lambda)^\ast$ with lowest weight $-\lambda$ and lowest weight vector $v_\lambda^\ast, $ satisfying the relations:
$$\lie n_{aff}^-.v_\lambda^\ast =0, \hspace{.75cm} h.v_\lambda^\ast = -\lambda(h)v_{\lambda}^\ast, \ \forall\ h\in\lie h_{aff}, \hspace{.75cm} (x_{\alpha_i}^+)^{\lambda(\alpha_i^\vee)+1}.v_\lambda^\ast =0, \ \forall\ 1\leq i\leq n+1.$$

We now record some results on integrable $\lie g_{aff}$-modules with 
finite-dimensional weight spaces that we shall need. Part(i) of the proposition
 was proved in \cite{C1}, part(ii) was proved in \cite{R3},
part(iii) was proved in \cite[Proposition 4.6]{DGK}, part(iv) was proved in \cite[Theorem 4.4]{A2} and part(v) was proved in \cite[Proposition 12.5(a)]{Kac}.

\begin{prop}Let $V$ be an integrable $\lie g_{aff}$-module having 
finite-dimensional weight spaces. 
\item[i.] If $V$ is an irreducible $\lie g_{aff}$-module on which the center 
acts by a positive (respectively negative) integer, then $V$ is isomorphic to $X(\Lambda)$ (respectively $X^\ast(\Lambda)$) for some $\Lambda\in P_{aff}^+$.
\item[ii.]If all eigenvalues of $K_1$ are non-zero, then $V$ is completely reducible as $\lie g_{aff}$-module.
\item[iii.] Given $\Lambda\in P_{aff}^+$, let $X^\omega(\Lambda)$ be the $\lie g_{aff}$-module whose underlying space is $X^\ast(\Lambda)$ and on which the action of $\lie g_{aff}$ is defined as : $$(x.f)(v) = -f(\omega(x).v), \hspace{2cm}\forall\ f\in X^\ast(\Lambda),\, x\in \lie g_{aff},\, v\in X(\Lambda),$$ where $\omega$ is the Cartan involution on $\lie g_{aff}$. Then,  $X(\Lambda)$ is isomorphic to $X^\omega(\Lambda)$ as $\lie g_{aff}$-modules.
\item[iv.] Given a finite-dimensional irreducible $\lie g_{fin}$-module $V(\lambda)$, let $\overline{V(\lambda)}:=V(\lambda)\otimes \C[t_1,t_1^{-1}]$,
be the $\lie g_{aff}$-module on which the action of $\lie g_{aff}$ is defined as follows: $$\begin{array}{ll}K_1. v\otimes t_1^s =0 &\ \forall\, v\in V(\lambda),\ s\in \Z,\\
d_1.v\otimes t_1^s = s v\otimes t_1^s &\ \forall\,  v\in V(\lambda),\ s\in \Z,\\ 
xt_1^r.v\otimes t_1^s = x.v \otimes t_1^{r+s} &\ \forall\, v\in V(\lambda),\ x\in \lie g_{fin}, \ s,r\in \Z. \end{array} $$ Then, for $\Lambda\in P_{aff}^+$, the $\lie g_{aff}$-module $X(\Lambda)\otimes \overline{V(\lambda)}$ is irreducible if 
\begin{equation}\lambda(\theta^\vee) > \Lambda(K_1)\hspace{.45cm} \text{and}\hspace{.45cm} \Lambda|_{\lie h_{fin}}(\theta^\vee)\leq  \Lambda(K_1) . \end{equation} 
\item[v.] Given $\lambda\in P_{aff}^+$ with $\lambda(K_1)\in \Z_{>0}$, $$P_{aff}(\lambda) = W_{aff}.\{\mu\in P_{aff}^+: \lambda\succeq \mu\},$$
where $W_{aff}$ denotes the Weyl group of $\lie g_{aff}$.\end{prop}

\noindent For $k>1$, $\lie g_{aff}$ will be identified the subalgebra $\lie g_{fin}\otimes \C[t_1,t_1^{-1}]\oplus \C K_1\oplus \C d_1$ of $\cal T_k(\lie g)$.

\vspace{.25cm}

\subsection{} For $\l\in \Z_+$, $\lambda\in P_{fin}^+$, and $\theta$ the highest root of $\lie g_{fin}$, the reducility of the $\lie g_{aff}$-module $X(l\Lambda_{n+1}+\lambda)\otimes \overline{V(\theta)}$ plays an important role in obtaining the block decomposition of $\cal J_{int}^+.$ Using the representation theory of vertex operator algebras, Adamovic proved a set of results in this direction in \cite{A1, A2}. We recall them here and shall  use them in \secref{Proof.Prop}.

\label{VOA} Given a vertex operator algebra $\mathbb V$, Zhu constructed an associative algebra $\ba(\mathbb V)$ and with any
$\mathbb V$-module $M$, Frenkel and Zhu associated an $\ba(\mathbb V)$-bimodule $\ba(M)$. Using the fact that the generalized 
Verma module for $\lie g_{aff}$ with highest weight $l\Lambda_{n+1}$, namely $M(l\Lambda_{n+1})$, and its unique irreducible quotient $X(l\Lambda_{n+1})$  have the structure of a vertex operator algebra and for $\mu\in P_{fin}$, the Verma module $M(l\Lambda_{n+1}+\mu)$ is a module for $M(l\Lambda_0)$, Adamovi$\acute{\text{c}}$ proved a set of results that we state. Part(i) of the following proposition was proved in \cite[Theorem 3.2, Remark 3.2]{A2}.  Part(ii) of \propref{VOA} is derived from the proof of \cite[Proposition 1.2, Theorem 4.1]{A1},  part(iii) follows from \cite[Theorem 3.2]{A1} and \cite[Lemma 1.1, Remark 1.1]{A1} and part(iv) follows from \cite[Proposition 3.2]{A2}.

\begin{prop} For $l\in \mathbb N$, denote the vertex operator algebras  $M(l\Lambda_{n+1})$ and $X(l\Lambda_{n+1})$ by $M_l$ and $X_l$ respectively.
\begin{enumerate}
\item[i.] Let $\mu\in P_{fin}^+$ be such that $(\mu|\theta) \leq l$. 
Then, $$\Hom_{\lie g_{aff}}(\overline{V(\mu)}\otimes X(l\Lambda_{n+1}),X(l\Lambda_{n+1}+\mu))\neq 0.$$

\vspace{.15cm}

\item[ii.] For $\mu\in P_{fin}^+$, $X(l\Lambda_{n+1}+\mu)$ is a $X_l$-module if and only if $(\mu|\theta)\leq l.$

\item[iii.]  Suppose $X(l\Lambda_{n+1}+\mu_1)$, $X(l\Lambda_{n+1}+\mu_2)$ and 
$X(l\Lambda_{n+1}+\mu_3)$ are $X_l$-modules with $\mu_i\in P_{fin}^+$ for $i=1,2,3$. Then there exists a linear isomorphism between the vector spaces 
$\Hom_{\lie g_{aff}}(\overline{V(\mu_2)}\otimes X(l\Lambda_{n+1}+\mu_1),X(l\Lambda_{n+1}+\mu_3))$ and 
$\Hom_{\ba(X_l)}(\ba (X(l\Lambda_{n+1}+\mu_1))\otimes_{\ba(X_l)} V(\mu_2),V(\mu_3)).$
 
\vspace{.15cm}
\item[iv.] If $ X(l\Lambda_{n+1}+\lambda)$ and $X(l\Lambda_{n+1}+\mu)$ are $X_l$-modules then 
$\ba (X(l\Lambda_{n+1}+\lambda))\otimes_{\ba(X_l)} V(\mu)$ is isomorphic to  
$\ba (X(l\Lambda_{n+1}+\lambda))\otimes_{\cal U(\lie g_{fin})} V(\mu)$ which is further isomorphic to $$\frac{V(\lambda)\otimes V(\mu)}{J(\lambda,\mu)}, $$ where
$J(\lambda,\mu) = \cal U(\lie g_{fin})\{ v_{\lambda}\otimes (x_{\theta}^+)^{l+1-(\mu_1|\theta)}v_2\ : \ v_2\in V(\mu)\}$. 
\end{enumerate}  
\end{prop}

\vspace{.1cm}

\begin{rem} It follows from parts (ii), (iii) and (iv) of the above proposition that $\Hom_{\lie g_{aff}}(\overline{V(\mu_2)}\otimes X(l\Lambda_{n+1}+\mu_1),X(l\Lambda_{n+1}+\mu_3))$  is non-zero, whenever $\mu_1,\mu_2, \mu_3 \in P_{fin}^+$ are such that  $(\mu_3|\theta)\leq l$ and $V(\mu_3)$ is an irreducible component of $V(\mu_1)\otimes V(\mu_2)$. 
\end{rem}

\subsection{}\label{tensor} We now state a set of results  on the irreducible components of tensor product of irreducible finite-dimensional $\lie g_{fin}$-modules.  They have been listed in the article \cite{icm.ku}.

\begin{prop} Let $\lambda, \mu\in P_{fin}^+$.
\begin{enumerate} 
\item[i.] For any $w\in W_{fin}$, let $\widehat{\lambda+w\mu}$ denote the unique element in $P_{fin}^+$ in the $W_{fin}$ orbit of $\lambda+w\mu$. The irreducible module $V(\widehat{\lambda+w\mu})$  occurs in the decomposition of $V(\lambda)\otimes V(\mu)$.
\item[ii.] If $\lie g_{fin}$ is not of type $G_2$, then for every positive root $\beta$, V($\lambda+\mu-\beta$)  occurs in the decomposition of $V(\lambda)\otimes V(\mu)$, whenever $\lambda$ and $\mu$ are non-zero dominant integral weights of $\lie g_{fin}$.
\item[iii.] For $\lambda\in P_{fin}^+$, let $S_{\lambda} = 
\{1\le i\le n: \lambda(\alpha_i^\vee) =0\}$. If $\lie g_{fin}$ is of type $G_2$, 
then for a positive root $\beta$, $V(\lambda+\mu-\beta)$ occurs in the 
decomposition of $V(\lambda)\otimes V(\mu)$, whenever $\lambda$ and $\mu$ are 
non-zero dominant integral weights of $\lie g_{fin}$ such that 
$\lambda+\mu-\beta\in P_{fin}^+$ and $S_\lambda \cup S_\mu \subseteq 
\{ 1\le i\le n: \beta-\alpha_i \notin R_{fin}^+\cup\{0\}\}$.
\end{enumerate}
\end{prop}

\subsection{} Using \propref{VOA} and \propref{tensor} we prove the following 
result which is crucial for the main theorem.

\begin{prop}  Let $\lambda\in P_{fin}^+$ and let  $\mu = l\Lambda_{n+1}+\lambda \in P_{aff}^+$. When $\lie g_{fin}$ is simply-laced, $l\geq \theta(\theta^\vee)$, and $\lambda\equiv \omega_i\mod Q_{fin}$ for $i\in J_0$ (respectively $\lambda\equiv 0\mod Q_{fin}$ ),  there exists a sequence $\lambda_1, \lambda_2,\cdots,\lambda_r\in P_{fin}^+$ with $\lambda_j(\theta^\vee)\leq l$ such that $\lambda_1 = \lambda$, $\lambda_r = \omega_i$ (respectively $\lambda_r=\theta$)  and for each $1\le j\le r-1$,
either $$ \Hom_{\lie g_{aff}}(\lie g_{aff}\otimes X(l\Lambda_{n+1}+\lambda_{j+1}),X(l\Lambda_{n+1}+\lambda_j)) \ne 0$$ or $$ \Hom_{\lie g_{aff}}(\lie g_{aff}\otimes X(l\Lambda_{n+1}+\lambda_{j}),X(l\Lambda_{n+1}+\lambda_{j+1})) \ne 0. $$
\end{prop}
\proof
Since $\lie g_{fin}$ is a highest weight irreducible $\lie g_{fin}$-module with highest weight $\theta$, $$\lie g_{aff} = \overline{V(\theta)}\oplus\C K_1 \oplus\C d_1, $$  where $\overline{V(\theta)} = V(\theta)\otimes \C[t_1^{\pm1}]$. Thus, given $\mu_i= l\Lambda_{n+1}+\lambda_i \in P_{aff}^+$, $i=1,2$, there exists a non-zero homomorphism $\phi:\lie g_{aff}\otimes X(\mu_1) \rightarrow X(\mu_2)$ only if
\begin{equation} \Hom_{\lie g_{aff}}(\overline{V(\theta)}\otimes X(l\Lambda_{n+1}+\lambda_1),X(l\Lambda_{n+1}+\lambda_2)) \ne 0. 
\label{hom.mu.1.2}\end{equation} 
By \remref{VOA}, \eqref{hom.mu.1.2} holds only if $l\geq \theta(\theta^\vee)$ and $V(\lambda_2)$ is an irreducible component of $V(\theta)\otimes V(\lambda_1)$ . 

Define a function $ht^w: P_{fin}^+ \rightarrow \N$ by $$ht^w(\sum_{i=1}^{n} b_i\omega_i) = \sum_{i=1}^{n} a_i^\vee b_i,$$ where $a_i^\vee$ is the label of the $i^{th}$-node in Table \ref{table:affine}. Notice that 
for $\lambda\in P_{fin}^+$, $l\Lambda_{n+1}+\lambda +n\delta_1\in P_{aff}^+$ only if $ht^w(\lambda)\le l$ and, for $\Lambda\in P_{aff}^+$, $$\Lambda(\alpha_{n+1}^\vee) + ht^w(\Lambda|_{\lie h_{fin}}) = \Lambda(K_1).$$ We prove the proposition case by case by applying induction on $ht^w(\lambda)$. 

We show that, given $\Lambda = l\Lambda_{n+1} + \lambda + n\delta_1\in P_{aff}^+$, with $\lambda = \Lambda|_{\lie h_{fin}}$ and $l\ge \theta(\theta^\vee)$,
there exists a sequence $\eta_{1,\lambda}, \cdots, \eta_{s,\lambda}\in P_{fin}^+$ with
$$\eta_{1,\lambda} = \lambda,  \hspace{.37cm} \eta_{s,\lambda} = \theta \ \text{or}\ \omega_i, \ \text{for some }\ i\in J_0,$$  for all $1\le j\le s,$ either
$  \eta_{j,\lambda} = \eta_{j+1,\lambda}\pm (\theta-\beta)$ for some $\beta\in R_{fin}^+$, or $\eta_{j,\lambda} = \eta_{j+1,\lambda}\pm w(\theta)$ for some $w\in W_{fin}$,
and $$ht^w(\eta_{j,\lambda})  \le ht^w(\lambda)\le \lambda(\theta^\vee).$$
Then, it would follow from \propref{tensor}, that either $V(\eta_{j,\lambda})$ is an irreducible component of $V(\theta)\otimes V(\eta_{j+1,\lambda})$ or $V(\eta_{j+1,\lambda})$ is an irreducible component of $V(\theta)\otimes V(\eta_{j,\lambda})$  and a repeated application of the above method would give the desired sequence of dominant integral weights $\lambda_1, \lambda_2,\cdots,\lambda_r\in P_{fin}^+$. 

Let $\lie g_{fin}$ be of type $A_n, n\geq 1$. Then, $\theta(\theta^\vee) =2$ and  $J_0 =\{1,\cdots,n\}$. Thus there is nothing to prove if $ht^w(\Lambda|_{\lie h_{fin}})=1$. Assume $ht^w(\Lambda|_{\lie h_{fin}})=2$, that is,  
$\Lambda|_{\lie h_{fin}}\in \{\omega_i+\omega_j: 1\leq i,j \leq n\}$. Then, using 
the relations, 
\begin{equation}\label{A.n} \begin{array}{c}
\omega_i+\omega_j +\alpha_{i+1}+\cdots+\alpha_{j-1} = \omega_{i-1}+\omega_{j+1}, \hspace{.25cm} \text{for}\ i\leq j ,\\
\omega_i+\alpha_1+\alpha_2+\cdots+\alpha_{i-1} = \omega_1+\omega_{i-1},\\
 \omega_k+\alpha_{k+1}+\cdots+\alpha_{n} = \omega_{k+1}+\omega_n,  \\
\omega_1+\omega_n = \theta,
 \end{array}
\end{equation}
we obtain the sequence $\{\eta_{r,\omega_i+\omega_j}\}_r$ as follows.\\
If $1\leq i,j\leq n$ is such that $i+j \leq n+1$ and $i\leq j$,  we have,
\begin{eqnarray*}
\eta_{r, \omega_i+\omega_j} = \omega_{i-r+1}+\omega_{j+r-1}, & \text{for} \
 1\le r\le i, \hspace{.45cm}
\eta_{i+1, \omega_i+\omega_j} = \left\{ \begin{array}{cc} \omega_{i+j}, & \text{if} \ i+j<n+1,\\
\theta, & \text{if} \ i+j=n+1 \end{array}
\right.\end{eqnarray*}
If $1\leq i,j\leq n$ is such that $i+j-(n+1) >0$ and $i\leq j$,  we have,
\begin{eqnarray*}
\eta_{r, \omega_i+\omega_j} = \omega_{i-r+1}+\omega_{j+r-1}, & \text{for} \
 1\le r\le n-j+1, \hspace{.45cm}
\eta_{n-j+2, \omega_i+\omega_j} = \omega_{i+j-(n+1)}.\end{eqnarray*}
Observe that $\eta_{r,\omega_i+\omega_j}(\theta^\vee)\leq 2$ for each $r$. Hence the sequence obtained satisfies the desired conditions. Now, if $ht^w(\Lambda|_{\lie h_{fin}})>2$, then there exists $\lambda', \lambda''\in P_{fin}^+$ such that $\Lambda|_{\lie h_{fin}} = \lambda'+\lambda''$ and $ht^w(\lambda'') =2$. Then, applying the above relations on $\lambda''$, we can obtain a sequence $\{\eta_{r,\Lambda|_{\lie h_{fin}}}\}_r$, such that for some $r\in \N$,
$$\eta_{r,\Lambda|_{\lie h_{fin}}} = \left\{ \begin{array}{cc}
\lambda'+\omega_i, & \text{if} \ \lambda''\equiv \omega_i\mod Q_{fin},\\ 
\lambda', & \text{if} \ \lambda''\equiv 0 \mod Q_{fin}
\end{array}\right. $$ 
Since $ht^w(\lambda'+\omega_i)< ht^w(\lambda)-1$ and $ht^w(\lambda')< ht^w(\lambda)-2$, applying induction and using the above relations repeatedly, we obtain the desired sequence.

Likewise, to prove the result in each case, it is sufficient to obtain the sequences $\{\eta_{r,\mu}\}_r$ in the case when $\mu$ is a fundamental weight of $\lie g_{fin}$ or  $\mu\in P_{fin}^+$ with $ht^w(\mu)\le 2$.  We thus list the required relations and case-wise, give the sequence $\{\eta_{r,\mu}\}_r$ in the cases when $\mu$ is a fundamental weight of $\lie g_{fin}$ or $\mu\in P_{fin}^+$ with $ht^w(\mu)\le 2$.

\vspace{.1cm}

For $\lie g_{fin}$ of type $D_n, n\geq 4$,  $\theta(\theta^\vee) =2$, $J_0 =\{1,n-1,n\}$, the set of $\lambda\in P_{fin}^+$ with  $ht^w(\lambda) =1$ is 
$\{\omega_1,\omega_{n-1}, \omega_n\}$ and the set of $\lambda\in P_{fin}^+$ with 
$ht^w(\lambda) =2$ is $ \{\omega_i: 2\leq i\leq n-2\}\cup \{\omega_i+\omega_j: i,j \in \{1,n,n-1\} \}$. Clearly there is nothing to prove when $ht^w(\lambda)=1$. When $ht^w(\lambda) =2$, using the relations,
\begin{equation}\label{D.n} 
\begin{array}{c} 
\omega_i = \omega_{i-2}+\alpha_{i-1} + 2(\alpha_i+\cdots+\alpha_{n-1})+\alpha_{n-1}+\alpha_n, \hspace{3\leq i\leq n-2},\\
\omega_1+\omega_n = \omega_{n-1}+\theta-(\underset{2\leq j\leq n-1}{\sum}\alpha_j), \hspace{.37cm}
\omega_1+\omega_{n-1} = \omega_n+ \theta -(\underset{2\leq j\leq n-2}{\sum}\alpha_j + \alpha_n)\\ 
 2\omega_1 = \omega_2+\alpha_1, \hspace{.37cm} 2\omega_{n-1}=\omega_{n-2}+\alpha_{n-1},\hspace{.37cm} 2\omega_{n}=\omega_{n-2}+\alpha_{n}\\
\omega_{n-1}+\omega_{n} = w_{n-3}+\alpha_{n-2}+\alpha_{n-1}+\alpha_{n}\hspace{.37cm} \omega_2=\theta,
 \end{array}
\end{equation}
for $2\le i\le n-2$, we get,
\begin{eqnarray*}\eta_{r,\omega_i} = \omega_{i-2(r-1)}, & \ \text{for}\ 1\le r\le i/2, \hspace{.15cm} \text{when $i$ is even},\\
\eta_{r,\omega_i} = \omega_{i-2(r-1)}, & \text{for}\ 1\le r\le (i+1)/2, \hspace{.15cm}, \hspace{.15cm} \text{when $i$ is odd}.
\end{eqnarray*}
Further we have,
\begin{eqnarray*}
\eta_{1,2\omega_1} = 2\omega_1, \eta_{2,2\omega_1} = \omega_2,\ & \eta_{1,2\omega_n} = 2\omega_n, \eta_{2,2\omega_n} = \omega_{n-1},\ &
 \eta_{1,2\omega_{n-1}} = 2\omega_{n-1}, \eta_{2,2\omega_{n-1}} = \omega_n,\end{eqnarray*}
\begin{eqnarray*}
\eta_{1,\omega_1 +\omega_n} = \omega_1 + \omega_n,\ \eta_{2,\omega_1+\omega_n} = \omega_{n-1},\ & \eta_{1,\omega_1+\omega_{n-1}} = \omega_1+ \omega_{n-1},\ 
\eta_{2,\omega_1+\omega_{n-1}} = \omega_n, \end{eqnarray*}
\begin{eqnarray*}
\eta_{1,\omega_{n-1}+\omega_n} = \omega_{n-1}+\omega_n, \hspace{.35cm}  \eta_{2,\omega_{n-1}+\omega_n} = \omega_{n-3},\\
\eta_{2+r,\omega_{n-1}+\omega_n} = \omega_{n-3-2r} , \left\{\begin{array}{ll} \text{for}\ 1\le r\le (n-4)/2,\ \text{when $n$ is even},\\
 \text{for}\ 1\le r\le (n-5)/2,\ \text{when $n$ is odd}.\\
\end{array}\right.
\end{eqnarray*}

\hspace{.1cm}

For $\lie g_{fin}$ of type $E_6$, $\theta(\theta^\vee)=2$ and $J_0=\{1,6\}$.  The set of $\lambda\in P_{fin}^+$ of $ht^w(\lambda) =1$ is $\{\omega_1, \omega_6\}$, the set of $\lambda\in P_{fin}^+$ of $ht^w(\lambda) =2$ is $\{\omega_2, \omega_3,\omega_5, \omega_1+\omega_6, 2\omega_1, 2\omega_6\}$ and the set of other fundamental weights and the corresponding value of $ht^w$ are, $\{\omega_4\}$ with 
$ht^w(\omega_4)= 3.$ 
Notice that there is nothing to prove in the case when $ht^w(\lambda) =1$. For $ht^w(\lambda)>1$, using the following relations,  
\begin{equation}\begin{array}{l} \label{E.6}
\omega_2=\theta, \hspace{.37cm} 2\omega_1=\omega_3+\alpha_1, \hspace{.37cm} 2\omega_6=\omega_5+\alpha_6 , \\
\omega_3=\omega_6+\alpha_1+\alpha_2+2\alpha_3+2\alpha_4+\alpha_5,\\
\omega_4=\omega_2+\alpha_1+\alpha_2+2\alpha_3+3\alpha_4+2\alpha_5+\alpha_6,\\
\omega_5=\omega_1+\alpha_2+\alpha_3+2\alpha_4+2\alpha_5+\alpha_6,\\
\omega_1+\omega_6 =\omega_2+\alpha_1+\alpha_3+\alpha_4+\alpha_5+\alpha_6
\end{array}\end{equation} 
and the sequences,
\begin{eqnarray*}
\eta_{1,\omega_3} = \omega_3, \hspace{.37cm} \eta_{2,\omega_3} = \omega_6,
\hspace{.9cm} \\
\eta_{1,\omega_4} = \omega_4, \hspace{.37cm} \eta_{2,\omega_4} = \omega_2, 
\hspace{.9cm}\\
\eta_{1,\omega_5} = \omega_5, \hspace{.37cm} \eta_{2,\omega_5} = \omega_1, 
\hspace{.9cm}\\ 
\eta_{1,\omega_1+\omega_6} = \omega_1+ \omega_6, \hspace{.37cm} 
\eta_{2,\omega_1+\omega_6} = \omega_2,\\
\hspace{.15cm}\eta_{1,2\omega_1} = 2\omega_1, \hspace{.37cm} \eta_{2,2\omega_1} = \omega_3,\hspace{.37cm} \eta_{3,2\omega_1} = \omega_6, \\ \hspace{.5cm}
\eta_{1,2\omega_6} = 2\omega_6, \hspace{.37cm} \eta_{2,2\omega_6} = \omega_5, \hspace{.37cm} \eta_{3,2\omega_6} = \omega_1,\end{eqnarray*}
one can obtain a sequence $\{\eta_{r,\lambda}\}_r$ of the desired form for all $\lambda\in P_{fin}^+$.

For $\lie g_{fin}$ of type $E_7$, $\theta(\theta^\vee)=2$ and $J_0=\{7\}$.  The set of $\lambda\in P_{fin}^+$ of $ht^w(\lambda) =1$ is $\{\omega_7\}$, the set of $\lambda\in P_{fin}^+$ of $ht^w(\lambda) =2$ is $\{\omega_1, \omega_2,\omega_6, 2\omega_7\}$ and the set of other fundamental weights and the corresponding value of $ht^w$ are, $\{\omega_3, \omega_5, \omega_4\}$ with $ht^w(\omega_3) = 3= ht^w(\omega_5)$ and $ht^w(\omega_4) =4$.  
As above, there is nothing to prove in the case when $ht^w(\lambda) =1$. For $\lambda\in P_{fin}^+$ with $ht^w(\lambda)>1$, using the following relations,
\begin{equation}\begin{array}{l} \label{E.7}
\omega_1=\theta, \hspace{.37cm} \omega_2=\omega_7+\alpha_1+2\alpha_2+2\alpha_3+3\alpha_4+2\alpha_5+2\alpha_6,\\
\omega_3=\omega_6+\alpha_1+\alpha_2+2\alpha_3+2\alpha_4+\alpha_5,\\
\omega_4=\omega_6+ [\theta -(\alpha_2+\alpha_3+2\alpha_4+\alpha_5)]\\
\omega_5 = \omega_2+\alpha_1 +\alpha_2+2\alpha_3+3\alpha_4+3\alpha_5+\alpha_6+\alpha_7,\\
\omega_6=\omega_1+\alpha_2+\alpha_3+2\alpha_4+2\alpha_5+2\alpha_6+\alpha_7,\\
2\omega_7 =\omega_6+\alpha_7,  
\end{array}\end{equation} 
and the sequences
\begin{eqnarray*}
\eta_{1,\omega_2} = \omega_2, \hspace{.37cm} \eta_{2,\omega_2} = \omega_7,
\hspace{.9cm} \\
\eta_{1,\omega_3} = \omega_3, \hspace{.37cm} \eta_{2,\omega_3} = \omega_6, \hspace{.37cm} \eta_{3,\omega_3} = \omega_1,\\
\eta_{1,\omega_4} = \omega_4, \hspace{.37cm} \eta_{2,\omega_4} = \omega_6, \hspace{.37cm} \eta_{3,\omega_4} = \omega_1,\\
\eta_{1,\omega_5} = \omega_5, \hspace{.37cm} \eta_{2,\omega_5} = \omega_2, \hspace{.37cm} \eta_{3,\omega_5} = \omega_7,\\ 
\eta_{1,\omega_6} = \omega_6, \hspace{.37cm} \eta_{2,\omega_6} = \omega_1,\hspace{.9cm}\\
\eta_{1,2\omega_7} = 2\omega_7, \hspace{.37cm} \eta_{2,2\omega_7} = \omega_6, \hspace{.37cm} \eta_{3,2\omega_7} = \omega_1,\end{eqnarray*}
one can obtain a sequence $\{\eta_{r,\lambda}\}_r$ of the desired form for all $\lambda\in P_{fin}^+$.

For $\lie g_{fin}$ of type $E_8$, $\theta(\theta^\vee) =2$ and $J_0=\emptyset$. In this case there does not exist $\lambda\in P_{fin}^+$ with $ht^w(\lambda) =1$. The set of $\lambda\in P_{fin}^+$ with $ht^w(\lambda)=2$ is $\{\omega_1,\omega_8\}$, and the set of other fundamental weights and the corresponding values for $ht^w$are $\{\omega_2,\omega_3, \omega_4,\omega_5, \omega_6,\omega_7\}$ with $ht^w(\omega_2) =3 =ht^w(\omega_7)$, $ht^w(\omega_3) =4 =ht^w(\omega_6)$, $ht^w(\omega_5) = 5$ and $ht^w(\omega_4) =6$.   Then using the following relations, 
\begin{equation}\begin{array}{l}
 \omega_1=\omega_8 + 2\alpha_1+2\alpha_2+3\alpha_3+4\alpha_4+3\alpha_5+2\alpha_6+\alpha_7,\\
\omega_2 =\omega_1 +\alpha_1 +3\alpha_2 +3\alpha_3+5\alpha_4+4\alpha_5+3\alpha_6+2\alpha_7+\alpha_8,\\
\omega_3=\omega_2+2\alpha_1+2\alpha_2+4\alpha_3+5\alpha_4+4\alpha_5+3\alpha_6+2\alpha_7+\alpha_8,\\
\omega_4=\omega_5+2\alpha_1+3\alpha_2+4\alpha_3+6\alpha_4+4\alpha_5+3\alpha_6+2\alpha_7+\alpha_8\\
\omega_5=\omega_6+2\alpha_1+3\alpha_2+4\alpha_3+6\alpha_4+5\alpha_5+3\alpha_6+2\alpha_7+\alpha_8\\
\omega_6=\omega_7+2\alpha_1+3\alpha_2+4\alpha_3+6\alpha_4+5\alpha_5+4\alpha_6+2\alpha_7+\alpha_8,\\
\omega_7 = \omega_8+2\alpha_1+3\alpha_2+4\alpha_3+6\alpha_4+5\alpha_5+4\alpha_6+3\alpha_7+\alpha_8,\\
\omega_8 =\theta, 
\end{array}\end{equation}  and the sequences
\begin{eqnarray*}
\eta_{1,\omega_1} = \omega_1, \hspace{.37cm} \eta_{2,\omega_1} = \omega_8,
\hspace{.9cm} \\
\eta_{1,\omega_2} = \omega_2, \hspace{.37cm} \eta_{2,\omega_2} = \omega_1, \hspace{.37cm} \eta_{3,\omega_2} = \omega_8,\\
\eta_{1,\omega_3} = \omega_3, \hspace{.37cm} \eta_{2,\omega_3} = \omega_2, \hspace{.37cm} \eta_{3,\omega_3} = \omega_1,\hspace{.37cm} \eta_{4,\omega_3} = \omega_8,\\
\eta_{1,\omega_4} = \omega_4, \hspace{.37cm} \eta_{2,\omega_4} = \omega_5, \hspace{.37cm} \eta_{3,\omega_4} = \omega_6,\hspace{.37cm} \eta_{4,\omega_4} = \omega_7, \hspace{.37cm} \eta_{5,\omega_4} = \omega_8,    \\ 
\eta_{1,\omega_5} = \omega_5, \hspace{.37cm} \eta_{2,\omega_5} = \omega_6,\hspace{.37cm} \eta_{3,\omega_5} = \omega_7,\hspace{.37cm} \eta_{4,\omega_5} = \omega_8,\\
\eta_{1,\omega_6} = \omega_6, \hspace{.37cm} \eta_{2,\omega_6} = \omega_7, \hspace{.37cm} \eta_{3,\omega_6} = \omega_8,\\
\eta_{1,\omega_7} = \omega_7, \hspace{.37cm} \eta_{2,\omega_7} = \omega_8,
\hspace{.9cm} 
\end{eqnarray*}
one can obtain a sequence $\{\eta_{r,\lambda}\}_r$ of the desired form for all $\lambda\in P_{fin}^+$.
\endproof

\label{algo}

\vspace{.15cm}

\section{The categories $\cal I$ and $\cal I_\fn$ }
\label{I.fin.ast}

\subsection{}  A  $\cal T_k(\lie g)$-module is said to be integrable if it is
$\lie h_\tr$-diagonalizable and the root vectors corresponding to the real roots of $\cal T_k(\lie g)$ are locally nilpotent on $V$. Thus an integrable $\cal T_k(\lie g)$-module is of the form $$V = \underset{\lambda\in \lie h_\tr^\ast}{\oplus} V_\lambda, \hspace{.5cm} \text{where}\ V_\lambda=\{v\in V: hv=\lambda(h)v,\ \text{for all}\ h\in \lie h_\tr\}. $$ We set $P(V) =\{\lambda\in \lie h_\tr^\ast: V_\lambda\neq 0\}$ as the set of all weights of an integrable $\cal T_k(\lie g)$-module $V$.
\label{integrable}

 For  $\beta \in R^{re}$ 
let $r_\beta(\lambda)= \lambda - \lambda(\beta^\vee)\beta,$ for $\lambda\in \lie h_\tr^\ast$. 
 Let $W_\tr= \langle r_\beta: \beta\in R_\tr^{re} \rangle, $ be the group generated by the operators $r_\beta$ for $\beta\in R_\tr^{re}.$  
Parts (i-ii) of the following is standard and part (iii) has been  proved in \cite{Kh}. (Refer to \cite{R2, Kh} for details.)

\label{sl2rep}
\begin{lem}Let $V$ be an integrable $\cal T_k(\lie g)$-module. 
     For all $\lambda\in P(V)$, the following hold. \begin{enumerate}
\item[i.]$\lambda(\alpha_i^\vee)\in \Z$, for $1\leq i\leq n+k.$
\item[ii.] $w\lambda \in P(V)$ and $\dim V_\lambda = \dim V_{w\lambda}$ for all $w\in W_\tr$.
\item[iii.] Let $\alpha\in R_\fn^+$ and $\beta=\alpha+m_i\delta_i \in R_\tr^{re}.$ Then, $$r_\alpha r_\beta(\lambda) = \lambda + \frac{2}{(\alpha|\alpha)}( m_i\lambda(K_i)) \alpha - (\lambda(\alpha^\vee)+\frac{2}{(\alpha|\alpha)} m_i\lambda(K_i))\delta_i.$$ In particular, if $\lambda+\underset{i=1}{\overset{k}{\sum}} r_i\delta_i\in P(V)$ is such that $\lambda(K_1) = m$ and $ \lambda(K_j)=0$ for $j=2,\cdots,k$, then there exists $\bom=(m_2,\cdots,m_k)\in \Z^{k-1}$ with $0\leq m_i<|m|$ for $2\leq i\leq k$ such that $\lambda +r_1\delta_1 +\underset{i=2}{\overset{k}{\sum}}  m_i\delta_i \in P(V).$\end{enumerate}
   \end{lem}
\vspace{.15cm}

\subsection {} Let $\cal I$ be the category of integrable $\cal T_k(\lie g)$-modules and morphisms $$\hom_{\cal I}(V,V') = \hom_{\cal T_k(\lie g)}(V,V'), \hspace{.35cm} V, V'\in \Ob \cal I.$$
For $\boa=(a_1,\cdots,a_k)\in \C^k$, set
$$V\{\boa\}=\{v\in V: d_iv=(a_i+r_i)v\ \text{for some}\ r_i\in \Z, \ 1\leq i\leq k\}.$$ Clearly, $V\{\boa\}$ is a $\cal{T}_k(\lie g)$-submodule of $V$ and $V\{\boa\}=V\{\bob\}$ if and only if $\boa-\bob\in \Z^k.$ For any $\bar{\boa}\in\C^k/\Z^k$, let $\cal I\{\bar{\boa}\}$ be the full subcategory of integrable $\cal T_k(\lie g)$-modules $V$ satisfying $V=V\{\boa\}$, where $\boa$ is any representative of $\bar{\boa}$. 
In \cite[Lemma 3.2]{CG1}, the following result was proved for graded level zero  integrable representations of affine Lie algebras. The proof for integrable representations of toroidal Lie algebras is analogous.

\begin{lem} Let ${\boa}, \bob\in \mathbb{C}^k$ be such that $\boa-\bob\not\in \Z^k$. Then for $V$ in $\cal I\{\bar\boa\}$ and $V^\prime$ in $\cal I\{\bar\bob\}$, $\Ext^1_{\cal I}(V,V^\prime)=0.$ In particular,
$$\cal I= \underset{\bar\boa\in\C^k/\Z^k}{\bigoplus} \cal I\{\bar\boa\}$$ and the categories $\cal I\{\bar\boa\}$ are equivalent for all $\bar\boa\in\C^k/\Z^k$.
\end{lem}

\noindent Without loss of generality we thus restrict ourselves to the subcategory $\cal I\{\bar{\bold{0}}\}$ of $\cal I$.

\subsection{} For $\bom\in \Z^k$ let  $\cal I^{(\bom)}$ be the full subcategory of $\cal I$ whose objects are $\cal T_k(\lie g)$-modules on which the zero degree central element $K_i$  acts by the integer $m_i$ for $1\leq i\leq k$ . Note that for a fixed $\bom\in \mathbb{Z}^k$, using the structure theory of finitely generated $\mathbb {Z}$-modules, one can find a basis $\{\boh_1,\boh_2,\cdots,\boh_k\}$ of $\mathbb{Z}^k$ such that $\bom=(m_1,m_2,\cdots,m_k) = m\boh_1$, where $m=\gcd(m_1,\cdots,m_k)$.  Hence, the following result. 

\begin{prop} \label{central.action}
Let $V$ be an integrable $\cal T_k(\lie g)$-module. Then $$V= \underset{\bom\in \Z^k}{\bigoplus} V^{(\bom)},$$ where 
  $V^{(\bom)}$ is an object of $\cal I^{(\bom)}$ for each $\bom\in \Z^k$. 
Furthermore,    $\Ext_{\cal I}^1(V,U)=0$ for all $V\in \Ob \cal I^{(\bom)}$ and $U\in \Ob \cal I^{(\bon)}$ whenever $\bom,\bon\in \Z^k$ are such that $\bom\neq \bon$. In particular,
   $$\cal I = \underset{\bom\in \Z^k}{\bigoplus} \cal I^{(\bom)},$$ and, for a fixed  $\bom=(m_1,\cdots,m_k)\in \Z^k-\{\bold{0}\}$ with $m = \gcd(m_1,\cdots,m_k)$,   the category $\cal I^{(\bom)}$ is equivalent to $\cal I^{(m\boe_1)}$ where $\{\boe_i: 1\le i\le k\}$ is the canonical basis of $\Z^k$ with $e_i = (0,\cdots,\underset{i^{th}}{1},0,\cdots,0)$. 
  \end{prop}

\subsection{} Let $\cal I_\fn^{(\bom)}$ be the full subcategory of $\cal I^{(\bom)}$ whose objects are $\cal T_k(\lie g)$-modules having finite-dimensional weight spaces. Let  \begin{equation} \cal I_\fn^\bast := \underset{\bom\in \Z^k-\bold{0}}{\bigoplus} \cal I^{(\bom)}_\fn. \end{equation} It follows from above that
$$ \cal I_\fn = \cal I_\fn^{(0)}\,  \bigoplus\, \cal I_\fn^\bast.$$
In the case when $k=1$, the structure of the category $\cal I_\fn^{(\bold{0})}$ was studied in \cite{CG1}. The same methods can be extended to study the structure of the category of $\cal T_k(\lie g)$-modules in $\cal I_{\fn}^{(\bold{0})}$ for any integer $k$. In view of  \propref{central.action}, we shall now try to extend these methods to study the structure of the category $\cal I_\fn^{(m\boe_1)}$.

\subsection{} The following result was proved in \cite[Proposition 2.4]{R2}, in the case when $V$ is an irreducible $\cal T_k(\lie g)$-module in $\cal I_{fin}^{(m\boe_1)}, m>0$. We now prove it for an arbitrary $\cal T_k(\lie g)$-module $V$ in $\cal I_{fin}^{(m\boe_1)}$.
\label{highest.wt}
\begin{prop}Let  $V$ be an  integrable $\cal T_k(\lie g)$-module in  $\cal I_{fin}^{(m\boe_1)}$, where $m>0$. Then  given $\mu\in P(V)$ with $\mu|_{\lie h_{aff}}\in P_{aff}^+$, there exists $\eta\in Q_{aff}^+$ such that $\mu+\eta\in P(V)$ but $\mu+\eta+\eta'\notin P(V)$ for all $\eta'\in Q_{aff}^+$.
\end{prop}
\proof For a contradiction suppose that there exists $\mu\in P(V)$ with $\mu|_{\lie h_{aff}}\in P_{aff}^+$, such that for each $\eta\in Q_{aff}^+$ satisfying $\mu+\eta\in P(V)$ there exists $\eta'\in Q_{aff}^+$ such that $\mu+\eta+\eta'\in P(V)$. Then there exists an infinite sequence $(\eta_r)_{r\geq 1}$ in $Q_{aff}^+$
such that $\eta_{r+1} \succ \eta_{r}$ and $\mu+\eta_r\in P(V)$ for all $r\geq 1.$ 
Set $W_r:=\cal U(\lie g_{aff}).V_{\mu+\eta_r}$. Since $V$ is an object of $\cal I_{fin}^{(m\boe_1)}$ and $W_r$ is a subspace of $V$, $W_r$ is an 
integrable $\lie g_{aff}$-module of non-zero level and by construction the set of $\lie h_{aff}$-weights of $W_r$ is a subset of $\mu|_{\lie h_{aff}}+\eta_r+Q_{aff}$. 
Hence, if $\lambda$ is an $\lie h_{aff}$-weight of $W_r$, 
then $\lambda + \sum_{i=2}^k (\mu|\Lambda_{n+i})\delta_i$ is a 
$\lie h_{tor}$-weight of $V$. This implies that the set of weight vectors of $W_r$ with $\lie h_{aff}$-weight $\lambda$ is a subset of the $\lie h_{tor}$-weight space $V_{\lambda +\sum_{i=2}^k(\mu|\Lambda_{n+i})\delta_i}$ of $V$. By assumption, $V$ has finite-dimensional weight spaces. Thus, it follows that $W_r$ has finite-dimensional $\lie h_{aff}$-weight spaces and, by \propref{affine}(ii), $W_r$ can be written as the direct sum of (possibly infinitely many) irreducible 
$\lie g_{aff}$-modules of the form $X(\mu_{r,s})$ with $\mu_{r,s}\in P_{aff}^+$. 
Choose $s_1$ such that $\nu_1:=\mu_{1,s_1}\succ \mu|_{\lie h_{aff}}$. Observe that 
such $s_1$ exists since by assumption $\mu|_{\lie h_{aff}}+\eta_1$ is a $\lie h_{aff}$ weight of $W_1$.  Now let $r_2$ be the 
smallest positive integer such that $\nu_2:=\mu_{r_2,s_2}\succ \mu|_{\lie h_{aff}}$ and  $\mu_{1,s_1}\neq \mu_{r_2,s_2}$. Note that such $r_2$  exists since
$\{\eta_r\}_{r\geq 1}$ is an increasing sequence in $Q_{aff}^+$ with respect to 
the ordering $\succeq$ and, as each weight space of $V$ is finite-dimensional, 
given any $r\geq 1$, by the above argument, there cannot exist infinity many direct
summands of $W_r$ having the same highest weight. Repeating the process, we 
obtain an infinite collection $\nu_r+\sum_{i=2}^k (\mu|\Lambda_{n+i})\delta_i \in P(V)$ such that $\nu_r\succ \mu|_{\lie h_{aff}}$, 
for $r\geq 1$ and there exists a $\lie g_{aff}$-submodule $W(\nu_r)$ of $V$ 
isomorphic to $X(\nu_r)$ as a $\lie g_{aff}$-module. By \propref{affine}(v), $\mu|_{\lie h_{aff}}$ is an $\lie h_{aff}$-weight of $W(\nu_r)$ for each $r\geq 1$. Since $\{\nu_r\}_{r\in \Z}$ is a set of  distinct $\lie h_{aff}$-weights and the sum of $W(\nu_r)$ is direct, this contradicts the finite-dimensionality of $V_\mu$.  Hence the proposition. 
\endproof

\subsection{}\label{V-aff+} Let  $\lie n_{aff}^+$ be the positive root space of the Lie subalgebra $\lie g_{aff} = \lie g\otimes \C[t_1^{\pm 1}]\oplus \C K_1\oplus \C d_1$ of $\cal T_k(\lie g)$. Given an integrable  $\cal T_k(\lie g)$-module $V$ define
\begin{equation} V^+_{aff}=\{v\in V_\lambda : \lie n^+_{aff}\otimes \C[t_2^{\pm 1},\cdots,t_k^{\pm 1}].v=0\},
\end{equation} i.e., $V_{aff}^+$ consists of all weight vectors $v\in V$ such that $x\otimes f.v=0$ for all $x\in \lie n_{aff}^+$ and $f\in \C[t_2^{\pm 1},\cdots,t_k^{\pm 1}].$
\vspace{.15cm}

\begin{cor} Let  $V$ be an integrable $\cal T_k(\lie g)$-module in  $\cal I_{fin}^{(m\boe_1)}$, $m>0$. Then, the set $V_{aff}^+$ is non-zero. Furthermore, 
there exists $\lambda\in \lie h_{tor}^\ast$ and $v\in V_{aff}^+$ such that $h.v = \lambda(h)v$ for all $h\in \lie h_\tr$, $\lambda|_{\lie h_{aff}}\in P_{aff}^+$ and $0\leq \lambda(d_j)< m$ for $2\leq j\leq k$.
\end{cor}
\proof It follows from the above proposition that given $\mu\in P(V)$ there exists $\eta_{\mu}\in Q_{aff}^+$ such that $\mu+\eta_\mu\in P(V)$ but $\mu+\eta_\mu+\eta \not\in P(V)$ for all $\eta\in Q_{aff}^+$. Now using the same proof as \cite[Proposition 2.4]{R2} and \lemref{sl2rep}(iii)  it will follow that there exists $\mu'\in P(V)$ with $(\mu'-\mu)|_{\lie h_{aff}}\in Q_{aff}^+$ and $0 \leq \mu'(d_j)<m$, $2\leq j\leq k$ such that $$\lie n_{aff}^+\otimes \C[t_2^{\pm1},\cdots,t_k^{\pm1}].V_{\mu'} =0,$$ which implies that $V_{\mu'}\cap V_{aff}^+$ is non-zero. \endproof

\subsection{} Using the analogues of \propref{highest.wt} and \corref{V-aff+}, the following result was proved in \cite{CG1} for an affine Lie algebra. It can be proved in the same manner for a $k$-toroidal Lie algebra.

\begin{thm} Let $V$ be an object in $\cal I_{fin}^{(0)}$. Then $V$ is isomorphic to a direct sum of indecomposable modules only finitely many of which are non-trivial.
\end{thm}

However it can be easily seen that the  $\cal T_k(\lie g)$-modules in 
$\cal I_{fin}^{(m\boe_1)}, m\neq 0$ do not in general satisfy the finite length 
property. For example, consider the $\cal T_k(\lie g)$-module 
$$V=\underset{s=0}{\overset{\infty}{\oplus}} X(\Lambda_i - s\delta_1)$$ when  
$\lie g=\lie{sl}_2(\C)$ and $i=0$ or 1. Clearly $V$ is an
integrable $\lie g_{aff}$-module with finite-dimensional weight spaces which
cannot be written as a direct sum of finitely many indecomposable
$\lie g_{aff}$-modules. As a first step towards understanding the category $\cal I_{fin}^{(m\boe_1)}$, we thus restricted our attention in this paper towards a subcategory $\cal J_{int}^{m}$  of $\cal I_{fin}^{(m\boe_1)}$.

\section{The Category $\cal J_{int}^\pm$} \label{J.fin}

\subsection{} Let $\cal J_{int}^{m}$ be the full subcategory of $\cal I_{int}^{(m\boe_1)}$ consisting of finite-length objects and let 
$$\cal J_{int}^{+} = \underset{m>0}{\bigoplus} \ \cal J^{m}_{int}.$$

 Clearly the simple objects of the category $\cal J_{int}^{m}$ are precisely the simple objects in $\cal I_\fn^{(m\boe_1)}$. Since they play an 
important role in determining  the structure of $\cal J_{int}^+$, we recall the results from 
\cite{R2,R05,Kh} that give a parametrization of the irreducibles in 
$\cal I_\fn^{\bast}$.

\label{irrep}

\subsection{} Let $\cal Z_1$ be the subspace of $\cal Z$ spanned by elements of
the form $t^{\bor}K_1$ with $\bor\in \Z^{k}_0$ (refer to \ref{section1.1} for the
notation) and let $\cal Z_1'$ be the subspace of $\cal Z$ spanned by the set of 
elements $\{t^\bor K_i, 2\leq i\leq k,\ \bor\in \Z^k\}\cup \{t^\bor K_1: \bor\in 
\Z^k-\Z^k_0\}$. Then $\cal Z= \cal Z_1 \oplus \cal Z_1'$ 
and  it is known from \cite[Theorem 4.5]{R2}, \cite[Proposition 4.1]{Kh} that, if $V$ is an irreducible
$\cal T_k(\lie g)$-module  in $\cal I_\fn^{(m\boe_1)}, m\neq 0$, then every element 
in $\cal Z_1'$ acts trivially on $V$. In other words every irreducible
$\cal T_k(\lie g)$-module  in $\cal I_\fn^{(m\boe_1)}, m\neq 0$, is in fact a 
module 
for the Lie algebra $\lie g_\fn \otimes \C[t_1^{\pm1},t_2^{\pm1},\cdots,t_k^{\pm1}]
\oplus \cal Z_1 \oplus D_1$.

\subsection{} Let $\cal L^c(\lie g) := \lie g_\fn \otimes \C[t_1^{\pm1},t_2^{\pm1},\cdots,t_k^{\pm1}] \oplus \cal Z_1 \oplus D_1 $. Then $\lie h_{aff}$ is an abelian subalgebra of $\cal L^c(\lie g)$. A $\cal L^c(\lie g)$-module $V$ is said to be integrable, if it is $\lie h_{aff}$-diagonalizable and  for all $\alpha\in R_{fin}$ and $\bom\in \Z^k$, the elements $x_{\alpha}\otimes t^\bom\in \cal L^c(\lie g)$, act locally nilpotently on $V$.

For an integrable $\cal L^c(\lie g)$ module $V$ set $$P^c(V) := \{\lambda\in P_{aff} : \text{there exists}\ v\in V\ \text{satisfying}\ h.v=\lambda(h)v, \ \text{for all}\ h\in \lie h_{aff} \}.$$

 \vspace{.15cm}\label{six.five}

 Let $\Pi$ be the monoid of finitely supported functions 
$\pi:\max \C[t_2^{\pm 1},\cdots,t_k^{\pm 1}] \rightarrow P_{aff}^+$, where, given 
$\pi,\pi'\in \Pi$ and $M\in \max \C[t_2^{\pm 1},\cdots,t_k^{\pm 1}]$, we define
$$(\pi+\pi')(M)= \pi(M)+\pi'(M),\hspace{.75cm} \supp(\pi)=
\{M\in \max \C[t_2^{\pm 1},\cdots,t_k^{\pm 1}]: \pi(M)\neq 0\},$$ and 
$$ \wt(\pi)=\sum_{M\in\supp(\pi)} \pi(M). $$
Given $\pi\in \Pi$ with $\supp(\pi) =\{M_1, M_2,\cdots, M_l\}$,  let 
$$X_\pi = X(\pi(M_1))\otimes\cdots\otimes X(\pi(M_l)). $$ One defines an
$\cal L^c(\lie g)$-module action on $X_\pi$ as follows:
$$Y\otimes f.  v_1\otimes \cdots \otimes v_l = \sum_{i=1}^l ev_{M_i}(f) v_1\otimes
\cdots\otimes Y.v_i\otimes\cdots\otimes v_l, $$ where $Y\in \lie g_{aff}$, 
$f\in \C[t_2^{\pm 1},\cdots,t_k^{\pm 1}]$, $v_i\in X(\pi(M_i))$ and  
$\ev_{M_i}: \C[t_2^{\pm 1},\cdots,t_k^{\pm 1}]\rightarrow \C$ is the evaluation map
at the point in $(\C^\ast)^{k-1}$ corresponding to the maximal ideal $M_i$ of 
$\C[t_2^{\pm 1},\cdots,t_k^{\pm 1}]$ for $1\leq i\leq l$. Since for each $i$, $X(\pi(M_i))$ is an integrable $\lie g_{aff}$-module, from the description of the $\cal L^c(\lie g)$-action on $X_\pi$, it is clear that $X_\pi$ is an integrable $\cal L^c(\lie g)$-module.  

Let 
$$ L(X_\pi) = X_\pi \otimes \C[t_2^{\pm 1},\cdots,t_k^{\pm 1}], $$
 be the $\cal T_k(\lie g)$-module on which the Lie algebra action is defined by :
 \begin{equation} \label{T.g.action.pi} 
\begin{array}{ll}Y\otimes f. (w\otimes f') = (Y\otimes f.w)\otimes ff', 
\hspace{1.5cm}& c.(w\otimes f') =0, \hspace{.35cm}  \forall\  c\in \cal Z_1',\\
 & \\
d_i. (w\otimes f') = w\otimes d_i(f'), \hspace{.35cm} \text{for}\ 2\leq i\leq k,
 \hspace{1.75cm}& d_1. w\otimes f' = d_1(w)\otimes f',\end{array}\end{equation}
where $Y\in \lie g_\fn\otimes \C[t_1^{\pm 1}] \oplus \C K_1,$ 
$f,f'\in \C[t_2^{\pm 1},\cdots,t_k^{\pm 1}]$ and $w\in X_\pi$. For 
$M\in \supp(\pi)$, let $v_{M}$ be the highest weight vector of $X(\pi(M))$ and 
let $$v_\pi:=  v_{M_1}\otimes v_{M_2}\otimes \cdots\otimes v_{M_l}.$$ Clearly 
$v_\pi\in L(X_\pi)_{aff}^+$ and 
$\cal U(\lie h_{aff}\otimes \C[t_2^{\pm 1},\cdots,t_k^{\pm 1}]).v_\pi$ is a 
$\Z^{k-1}$-graded subspace of $L(X_\pi)_{aff}^+$. Let 
$$\operatorname{Ann}_{\lie h_{aff}\otimes 
\C[t_2^{\pm 1},\cdots,t_k^{\pm 1}]} (v_\pi) = \{a \in \cal U(\lie h_{aff}\otimes 
\C[t_2^{\pm 1},\cdots,t_k^{\pm 1}]) |\ a.v_\pi =0\}$$
and  
$$\bold{A}_\pi := \cal U(\lie h_{aff}\otimes 
\C[t_2^{\pm 1},\cdots,t_k^{\pm 1}])/ \operatorname{Ann}_{\lie h_{aff}\otimes 
\C[t_2^{\pm 1},\cdots,t_k^{\pm 1}]} (v_\pi) .$$ 
It is easy to see that $\bold{A}_\pi$ 
is a $\Z^{k-1}$-graded algebra and  
$$\bold{A}_\pi= \underset{\bom\in \Z^{k}_0} {\oplus} 
\bold{A}_\pi[{\bom}], $$ where for each ${\bom}\in 
\Z^k_0$, $$\bold{A}_\pi[{\bom}]= 
\{h_1\otimes t^{{\bom^1}} h_2\otimes t^{{\bom^2}}\cdots \in 
\cal U(\lie h_{aff}\otimes \C[t_2^{\pm 1},\cdots,t_k^{\pm 1}])\ \text{with}\ \bom^i\in \Z^k_0,\ \text{and} \
\sum_{i}{\bom^i}= {\bom}\}.$$  
 For $\pi\in \Pi$ let  $$G_\pi:=
\{{\bom} =(0,m_2,\cdots,m_k)\in \Z^{k}_0: A_\pi[{\bom}] 
\neq 0\}.$$ 
From the definition of the action of $\cal T_k(\lie g)$ on $L(X_\pi)$ 
\eqref{T.g.action.pi} it is clear that \begin{equation}\label{T(h).action}
    h\otimes t^{\bom}.v_\pi = (\sum_{i=1}^r \pi(M_i)(h) \ev_{M_i}(t^{\bom})) 
v_\pi\otimes t^{\bom}=0, \hspace{.5cm}\forall\ {\bom}\in \Z_0^{k}-G_\pi.
    \end{equation} As $v_\pi\ne 0$,  it thus follows from  \eqref{T.g.action.pi} and  the definition of the $\cal L^c(\lie g)$-module action on $X_\pi$ that, $$ (\sum_{i=1}^r \pi(M_i)(h) \ev_{M_i}(t^{{\bom}})) =0, 
\hspace{.5cm}\forall\ {\bom}\in \Z^{k}-G_\pi.$$ It has been shown in \cite{Kh} 
that $G_\pi$ is a subgroup of $\Z^{k-1}$ of finite index. We shall refer to $G_\pi$ as the group associated with $\pi\in \Pi$ and denote the corresponding quotient group by $G^\pi = \Z^{k-1}/G_\pi$.
 
The following results have been proved in \cite[Proposition 3.5, Theorem 3.18, 
Example 4.2]{R2}, \cite{R05}.  

\begin{prop} For $\pi\in \Pi$, let $v_\pi$ be the highest weight vector of the $\cal L^c(\lie g)$-module $X_\pi$. Then the following hold.
\item[i.] Given $\bog\in \Z^k_0$  
$$X_\pi^{\bog} = \cal U(\cal T_k(\lie g)).v_\pi\otimes t^{\bog},$$  
is an irreducible $\cal T_k(\lie g)$-module.
\item[ii.] Every irreducible $\cal T_k(\lie g)$-module of non-zero level which has finite-dimensional weight space is isomorphic to $X_\pi^\bog$ for some $\pi\in \Pi$ and $\bog\in \Z_0^k$.
\item[iii.] $L(X_\pi)$ is completely reducible as a $\cal T_k(\lie g)$-module. In fact, as a $\cal T_k (\lie g)$-module, $L(X_\pi)$ is isomorphic to the direct sum of the irreducible $\cal T_k (\lie g)$-modules $X_\pi^\bog$, where the direct sum is taken over representatives of distinct cosets of $G_\pi$ in $\Z^{k-1}$. In other words, for $\bog, \bop \in\Z^k_0$ with $\bog-\bop\in G_\pi$, $X_\pi^\bog = X_\pi^\bop$ and $$L(X_\pi) ={\oplus} X_\pi^\bog,$$ where each summand 
$X_\pi^\bog$ appears with multiplicity one.     \end{prop}

\subsection{} Notice that for $\bob=(b_2,\cdots,b_k)\in (\C^\ast)^{k-1}$ the map 
$s_\bob:\C[t_2^{\pm 1},\cdots,t_k^{\pm 1}] \rightarrow 
\C[t_2^{\pm 1},\cdots,t_k^{\pm 1}]$ given by $t_i\mapsto b_it_i$, $i=2,\cdots,k$  
is an isomorphism. Denote by $\bob.M$ the image of a maximal ideal $M$ of 
$\C[t_2^{\pm 1},\cdots,t_k^{\pm 1}]$ under the isomorphism $s_\bob$. Define an 
action of $(\C^\ast)^{k-1}$ on $\Pi$ by:
$$\bob.\pi (M) = \pi(\bob.M), \hspace{1.5cm} \text{for all}\ M\in
\max\C[t_2^{\pm 1},\cdots,t_k^{\pm 1}].$$ With this notation we have the following
 result from \cite{Kh} on the isomorphism classes of irreducible 
representations of $\cal T_k(\lie g)$.

\label{isomorphism}

\begin{prop} Given $\pi,\pi'\in \Pi$, $\bog, \bog' \in \Z^k_0$, 
the irreducible $\cal T_k(\lie g)$-modules $X_\pi^{\bog}$ and $X_{\pi'}^{\bog'}$ are 
isomorphic if and only if there exists $\bob\in (\C^\ast)^{k-1}$ such that
\item[i.]$\supp(\pi')= \{\bob.M: M\in \supp(\pi)\}.$
\item[ii.] For each $M\in \supp(\pi)$, there exists one-dimensional 
$\lie g_{aff}$-module $Z_M$ such that  $X(\pi(M))\otimes Z_M$ is isomorphic to 
$X(\pi'(\bob.M))$ as a $\lie g_{aff}$-module. 
\item[iii.]$\bog-\bog'\in G_\pi$.
\end{prop}

\vspace{.1cm}

\noindent From the description of the action of $\cal L^c(\lie g)$ on $X_\pi$, 
it is clear that for all $\pi\in \Pi$, $X_\pi$ is a representation for the 
quotient $\underset{M\in \supp(\pi)}{\bigoplus} 
(\lie g_\fn \otimes \C[t_1^{\pm1},t_2^{\pm 1},\cdots,t_k^{\pm 1}]/M \oplus 
{\cal Z_1}) \oplus \C d_1$ of the Lie algebra $\cal L^c(\lie g)$ and  
the following is an immediate consequence of \propref{isomorphism}.
\vspace{.1cm}
 
\begin{cor}  Given $\pi,\pi'\in \Pi$, the irreducible 
$\cal L ^c(\lie g)$-modules $X_\pi$ and $X_{\pi'}$ are isomorphic if and only if 
$\pi'= \bob.\pi $ for some $\bob=(b_2,\cdots,b_k)\in (\C^\ast)^{k-1}$. i.e., 
there exists $\bob=(b_2,\cdots,b_k)\in (\C^\ast)^{k-1}$ such that
  $\supp(\pi_2)= \{\bob.M: M\in \supp(\pi_1)\},$ and for each 
$M\in \supp(\pi_1)$, there exists one-dimensional $\lie g_{aff}$-module $Z_M$ 
such that  $X(\pi_1(M))\otimes Z_M$ is isomorphic to $X(\pi_2(\bob.M))$ as a 
$\lie g_{aff}$-module.
 \end{cor}

\subsection{} For any $\cal T_k(\lie g)$-module $V$ in $\cal I_{fin}^{(m\boe_1)}$, let $V^\vee = \oplus_{\mu\in P(V)} V_{\mu}^\ast \subseteq V^\ast$ be the graded dual. Define a $\cal T_k(\lie g)$-module structure on $V^\vee$ by 
$$(x\otimes t^{\bom}.f)v= -f(\omega(x)\otimes t^{\bom}.v), \hspace{.45cm} x\in \lie g_{aff} ,\, \bom\in\Z^k_0, \, v\in V,$$ 
where $\omega$ is the Cartan involution on $\lie g_{aff}$ (\secref{affine}). With this $\cal T_k(\lie g)$-module structure we denote $V^\vee$ by $V^{\tilde{\omega}}.$ Clearly the $\cal T_k(\lie g)$-module $V^{\tilde{\omega}}$ is an object of $\cal I_{fin}^{(m\boe_1)}$. Using the finite-dimensionality of the weight spaces of objects in $\cal I_{fin}^{(m\boe_1)}$ it can be seen that if $V$ has finite-length then $V^{\tilde{\omega}}$ also has finite-length. Thus the functor sending $V$ to
 $V^{\tilde{\omega}}$ is exact and contravariant on $\cal J_{int}^+$. Furthermore using \propref{affine}(iii), it follows that for $\pi\in \Pi$, 
$L(X_\pi)$ is isomorphic to $L(X_\pi)^{\tilde{\omega}}$ as $\cal T_k(\lie g)$-module. Hence there exists $\pi^{\tilde{\omega}}\in \Pi$ such that $L(X_\pi)^{\tilde{\omega}}$ is isomorphic to $L(X_{\pi^{\tilde{\omega}}})$ and,  if $\pi_1,\pi_2\in \Pi$ are such that $\wt(\pi_1)-\wt(\pi_2)\in Q^+\backslash \{0\}$, then $\wt({\pi_2^{\tilde{\omega}}})-\wt(\pi_1^{\tilde{\omega}})\notin Q^+$.

\section{Highest $\ell$-weight integrable modules in $\cal I_\fn^\bast$ }

We now recall from \cite{RSF} results on highest weight 
indecomposable $\cal T_k(\lie g)$-modules in $\cal I_{fin}^\bast$.
\label{Weyl}
\subsection{}Given a Cartan subalgebra $\lie h_{aff}$ of $\lie g_{aff}$, set 
$$\begin{array}{c}\cal T_{k-1}(\lie h_{aff}):= \lie h_{aff}\otimes 
\C[t_2^{\pm1},\cdots,t_k^{\pm1}] \oplus \cal Z\oplus D_k\\ \\ 
\cal L^c(\lie h):=\lie h_\af\otimes \C[t_2^{\pm1},\cdots,t_k^{\pm1}] \oplus \cal 
Z_1\oplus \C d_1  \end{array}$$
\begin{defn} Let $V$ be an integrable $\cal T_k(\lie g)$-module in 
$\cal I_\fn^\bast$. A non-zero weight vector $v\in V_{aff}^+$ of weight $\Lambda$ 
is said to be a highest $\ell$-weight vector, if $\cal U(\cal T_{k-1}(\lie h_{aff})).v$ is 
an indecomposable $\cal T_{k-1}(\lie h_{aff})$-submodule  of $V$ and 
$\dim (\cal U(\cal T_{k-1}(\lie h_{aff}))v)_{\Lambda+\delta_{\bom}}\leq 1$ for all 
$\bom\in \Z^{k}_0.$ A $\cal T_k(\lie g)$-module is said to be 
highest $\ell$-weight if $V= \cal U(\cal T_k(\lie g))v$ for some highest 
$\ell$-weight vector $v\in V$.
  \end{defn}

Clearly,  $X_\pi^{\bog}$ is a highest $\ell$-weight modules for all $\pi\in \Pi$ and ${\bog}\in \Z^{k}_0$.

\subsection{}\label{weyl} In the category $\cal I_\fn^{(\bold{0})}$, the notion of Weyl
 modules was defined in \cite{CPweyl} in the case $k=1$. In \cite{RSF} the notion of Weyl modules was defined for
 a Lie algebra of the form $\lie G\otimes A$, where $\lie G$ is a Kac-Moody Lie
 algebra and $A$ is a commutative associative algebra with unity. We recall 
here the definition and properties of these modules. We shall need them in 
\secref{BLK}.

\begin{prop}For $\pi\in \Pi$ with $\wt(\pi) = \Lambda$, let $W_\pi$ be the $\cal L^c(\lie g)$-module 
generated by a vector $w_\pi$ satisfying the following conditions:
$$\begin{array}{ll} \lie n_{aff}^+\otimes \C[t_2^{\pm1},\cdots,t_k^{\pm1}].w_\pi =0,
 & h.w_\pi = \Lambda(h)w_\pi, \hspace{.15cm}\forall\ \ h\in \lie h_{aff} \\
(x_{\alpha_i}^-)^{\Lambda(\alpha_i)^\vee+1}.w_\pi =0, &\text{for}\ i=1,2\cdots,n+1\\ 
h\otimes t^{{\bom}}.w_\pi = (\sum_{M\in \supp(\pi)} \ev_M(t^{{\bom}})
\pi(M)(h))w_\pi, &\forall\ h\otimes t^{{\bom}}\in \lie h_{aff}\otimes 
\C[t_2^{\pm1},\cdots,t_k^{\pm1}].\end{array}$$
\item[(i)] $W_\pi$ is an integrable $\cal L^c(\lie g)$-module with 
finite-dimensional weight spaces.
\item[(ii)]$X_\pi$ is the unique irreducible quotient of $W_\pi.$
\item[(iii)] Let $V$ be an integrable $\cal L^c(\lie g)$-module with 
finite-dimensional weight spaces generated by a weight vector $v$ such 
that 
$$ \begin{array}{ll}\lie n_{aff}^+\otimes \C[t_2^{\pm1},\cdots,t_k^{\pm1}].v =0,
 & h.v = \Lambda(h)v, \hspace{.15cm}\forall\ \ h\in \lie h_{aff} \\
(x_{\alpha_i}^-)^{\Lambda(\alpha_i)^\vee+1}.v =0, &\text{for}\ i=1,2\cdots,n+1\\ 
h\otimes t^{{\bom}}.v = (\sum_{M\in \supp(\pi)} \ev_M(t^{{\bom}})
\pi(M)(h))v, &\forall\ h\otimes t^{{\bom}}\in \lie h_{aff}\otimes 
\C[t_2^{\pm1},\cdots,t_k^{\pm1}].\end{array}$$
Then $V$ is isomorphic to a quotient of $W_\pi.$
\item[(iv)] Suppose $\pi_1,\pi_2\in \Pi$ are such that $\pi = \pi_1 +\pi_2$ 
and $\supp(\pi_1)\cap \supp(\pi_2) = \emptyset$, then 
$$W_\pi \cong_{\cal L^c(\lie g)} W_{\pi_1}\otimes W_{\pi_2}.$$
\item[(v)] For $M\in \max \C[t_2^{\pm1},\cdots,t_k^{\pm1}]$, let $\pi_{\Lambda,M}\in \Pi$ be the function such that  $\supp(\pi_{\Lambda,M}) = \{M\}$ and $\wt(\pi_{\Lambda,M})= \Lambda.$ Then, $W_{\pi_{\Lambda, M}}$ is spanned by elements of the form 
$$ (Y_1\otimes f_1^{r_1})(Y_2\otimes f_2^{r_2})\cdots (Y_s\otimes f_s^{r_s})w_{\pi_{\Lambda,M}},$$ where $s\in \N$, $Y_i\in\lie g_{aff}$, $f_i\in M$ and, $0\le r_i <\Lambda(\beta_i^\vee)$ when $Y_i = x_{\beta_i}^-\otimes g_i$ with $\beta_i\in R_{fin}^+$ and $g_i\in \C[t_1^{\pm 1}]$ and $0\le r_i <\underset{1\le j\le n}{\max}\Lambda(\alpha_j^\vee)$ when $Y_i \in \lie h_{fin}\otimes \C[t_1^{\pm1}]$ for $i=1,\cdots,s$.
\end{prop}
\proof Part (i) is proved in \cite[Lemma 4.3]{Kh}. For part (ii), note that  
$w_\pi\in (W_\pi)_{aff}^+$, $P^c(W_\pi) \subseteq \Lambda -Q_{aff}^+$ and $\dim (W_\pi)_\Lambda =1$. Hence $W_\pi$ has a unique irreducible quotient. Since, by \propref{six.five} and \corref{isomorphism}, the orbits of $\Pi$ under the natural action of $\C^{k-1}$ on it,  parametrizes the isomorphism classes of irreducible integrable $\cal L^c(\lie g)$-modules with finite-dimensional weight spaces, it follows from \secref{six.five} that $X_\pi$ is the unique irreducible quotient of $W_\pi.$
\item[(iii)] Since $V$ is an integrable $\cal L^c(\lie g)$ and $\alpha_i$ for $1,2,\cdots,n+1$ are real roots, using the representation theory of $\lie {sl}_2(\C)$ it is easy to see that if $h.v= \Lambda(h)v$ for all $h\in \lie h_{aff}$ then $(x_{\alpha_i}^-)^{\Lambda(\alpha_i^\vee)+1}v =0 $ for $i=1,\cdots,n+1.$ Hence it follows that $V= \cal U(\cal L^c(\lie g)).v$ is a quotient of $W_\pi$.
\item[(iv)] This was proved in a more general setup in \cite[Section 4, Section 5]{RSF}.
\item[(v)] It was shown in \cite[Proposition 4.3]{Kh} that $W_\pi$ is spanned 
by elements of the form $$ (Y_1\otimes (t^{\bom_1})^{r_1})(Y_2\otimes (t^{\bom_2})^{r_2})\cdots (Y_s\otimes (t^{\bom_s})^{r_s})w_{\pi_{\Lambda,M}},$$ where $s\in \N$, $Y_i\in\lie g_{aff}$, $\bom_i\in \Z^k_0$ and, $0\le r_i <\Lambda(\beta_i^\vee)$ when $Y_i = x_{\beta_i}^-\otimes g_i$ with $\beta_i\in R_{fin}^+$ and $g_i\in \C[t_1^{\pm 1}]$ and $0\le r_i <\underset{1\le j\le n}{\max}\Lambda(\alpha_j^\vee)$ when $Y_i \in \lie h_{fin}\otimes \C[t_1^{\pm1}]$ for $i=1,\cdots,s$. Further by \cite[ Section 5]{RSF}, $$x_\beta^-\otimes f^{\Lambda(\beta)}.w_{\pi_{\Lambda,M}} =0, \hspace{.35cm} \forall \ f\in M,$$ and all positive real roots $\beta$ of $\lie g_{aff}$. These together imply the assertion of part(v) of the proposition.\endproof
 
\rem  From part(iv) of the above proposition, it is clear that, if $\supp(\pi)=\{M_1,\cdots,M_r\}$ and $W_\pi \cong_{\cal L^c(\lie g)} \underset{1\le j\le r}{\otimes} W_{\pi_j}$, with $\pi_j\in \Pi$ such that $\supp(\pi_j) = \{M_j\}$ and $\pi_j(M_j)=\pi(M_j)$ for $j=1,\cdots,r$, then $W_\pi$ is spanned by elements of the form $w_1\otimes w_2\otimes..\cdots\otimes w_r$ where $w_j\in W_{\pi_j}$ and, hence, $\wt(w_j)\subseteq \pi(M_j)-Q_{aff}^+$, for $j=1,\cdots,r.$

\begin{cor} Suppose $\psi\in \Pi$ is such that $X_\psi$ is a $\cal L^c(\lie g)$ irreducible constituent of $W_\pi$ for $\pi\in \Pi$. Then, $\supp (\psi) = \supp (\pi)$ and, for each $M\in \supp (\pi)$, $\pi(M)-\psi(M)\in Q_{aff}^+$. 
\end{cor}
\proof Follows from part(iv) and part(v) of the above proposition.\endproof

 \subsection{} Given a $\cal L^c(\lie g)$-module $V$, let $L(V) := V\otimes \C[t_2^{\pm1}, \cdots, t_k^{\pm1}]$ and define a $\cal T_k(\lie g)$-module structure on $L(V)$ as follows : 
 $$(Y\otimes f)(w\otimes a)= (Y\otimes f.w)\otimes fa, \hspace{.5cm} d_1(w\otimes f) = d_1(w)\otimes f, \hspace{.5cm} d_i(w\otimes f) = w\otimes d_i(f), \hspace{.25cm} \forall\ i=2,\cdots,k,$$ where $Y\otimes f\in \lie g_\fn \otimes \C[t_1^{\pm1},\cdots,t_k^{\pm1}]\oplus \cal Z_1$, $w\in  V$ and $f,a\in \C[t_2^{\pm1},\cdots,t_k^{\pm1}].$  

Observe that $L$ defines a functor from the category of integrable $\cal L^c(\lie g)$-modules to the category of integrable $\cal T_k(\lie g)$-modules and direct sums and short exact sequences are preserved under the functor 
$L$. 
  
\label{gr.weyl} 
 \begin{prop} Let $\pi\in \Pi$ and $L^{{\bog}}(W_\pi) = \cal U(\cal T_k(\lie g)).w_\pi\otimes t^{{\bog}}$, for  ${\bog} \in \Z^{k}_0$.
\begin{enumerate} 
\item[i.] For ${\bog},{\bog}'\in \Z^k_0$, $L^{{\bog}}(W_\pi) = L^{{\bog}'}(W_\pi)$ if and only if ${\bog}-{\bog}'\in G_\pi$. Furthermore, identifying $\Z_0^k$ with $\Z^{k-1}$, we see that 
  $$L(W_\pi) = \oplus L^{{\bog}}(W_\pi),$$ where the direct sum is taken over representatives of distinct cosets of $G_\pi$ in $\Z^{k-1}$. 
\item[ii.] Let $V$ be a highest $\ell$-weight $\cal T_k(\lie g)$-module generated by a vector $v\in V_{aff}^+$ such that
 \begin{equation}\label{L.W.pi}\begin{array}{c}
h.v = \wt(\pi)(h)v, \hspace{.15cm}\forall\ h\in \lie h_{aff}, \hspace{.35cm}  d_i(v) = g_i v, \hspace{.15cm} \ \text{for}\ i=2,3,\cdots,k\\ 
 h\otimes f.v = (\underset{M\in \supp(\pi)}{\sum} \ev_M(f)\pi(M)(h) )v\otimes f, \hspace{.35cm} \forall \ h\otimes f\in \lie h_{aff}\otimes \C[t_2^{\pm1},\cdots,t_k^{\pm1}].
\end{array}
\end{equation}   Then $V$ is a quotient of $L^{{\bog}'}(W_\pi)$ where ${\bog}'\in \Z^{k}_0$ is such that ${\bog}' -(0,g_2,\cdots,g_k)\in G_\pi.$ 
\item[iii.] If $V$ is a $\cal L^c(\lie g)$-module quotient of $W_\pi$ and $v$ is the image of $w_\pi$ in $V$, then for $\bos\in \Z^k_0$, $L^\bos(V) = \cal U(\cal T_k(\lie g))(v\otimes t^\bos)$ is a highest $\ell$-weight module and identifying $\Z_0^k$ with $\Z^{k-1}$, we have $$L(V) = \oplus L^\bos(V),$$ where the direct sum is taken over representatives of distinct cosets of $G_\pi$ in $\Z^{k-1}.$
\end{enumerate}
 \end{prop}
 
Considering the isotypical components of the finite abelian group $G^\pi$, the first  part of the proposition can be proved in the same way as  \cite[Proposition 5.5(ii)]{CG1}. For the second part, note that if $V$ is a highest $\ell$ weight module generated by a vector $v$ satisfying the conditions \eqref{L.W.pi}, then similar arguments as in \cite[Section 4.3]{Kh} show that one can uniquely associate with $V$ a quotient of the $\cal L^c(\lie g)$-module $W_\pi$. Now using the first part of the proposition, it will then follow that $V$ is quotient of $L^{{\bog}}(W_\pi)$ for some ${\bog}\in \Z^{k}_0.$ For part(iii) of the propsition observe that if $V$ is a $\cal L^c(\lie g)$-module quotient of $W_\pi$, then by the exactness of the functor $L$, $L(V)$ is a quotient of $L(W_\pi)$. Now part (iii) follows from part(i) of the proposition.

\section{Vanishing of $\Ext^1$ in $\cal J_{int}^{\pm}$}
\label{BLK}

\subsection{} Let $\Xi$ be the set of all finitely supported functions $\xi: \max \C[t_2^{\pm 1},\cdots,t_k^{\pm 1}] \rightarrow \Z\times \Gamma.$ Since $\Z$ and $\Gamma$ are abelian groups with respect to addition operation, regarding $\Z\times \Gamma$ as the direct product of abelian groups, it is easy to see that addition of functions defines on $\Xi$ the structure of an additive abelian group in an obvious way. Denoting the images of the fundamental weights $\{\omega_i\}_{1\leq i\leq n}$ in $\Gamma$ by $\{\bar{\omega}_i\}_{1\leq i\leq n}$, for for $M\in \max\C[t_2^{\pm 1},\cdots,t_k^{\pm 1}]$, define
$$\xi_{n+1,M}(S) = \left\lbrace \begin{array}{ll}
(1,0) & \text{if}\ S=M,\\
(0,0) & \text{otherwise} 
\end{array}\right. \hspace{3.65cm} $$
$$\xi_{i,M}(S) =\left\lbrace \begin{array}{ll}
(a_i^\vee,\bar{\omega_i}) & \text{if}\ S=M,\\
(0,0) & \text{otherwise}
\end{array}\right. \hspace{1.5cm} \ \forall \ 1\leq i \leq n, $$ 
where $a_i^\vee$ is the integer labelling the
ith-node in the Dynkin diagrams given in Table \ref{table:affine}.

Clearly $\Xi$ is a free abelian group generated by the set of elements $\{\xi_{i,M}: 1\leq i\leq n+1,M\in \max\C[t_2^{\pm 1},\cdots,t_k^{\pm 1}]\}$. Define an action of $(\C^\ast )^{k-1}$ on $\Xi$ by $$(\bob.\xi)(S) = \xi(\bob.S),$$ where $\bob=(b_2,\cdots,b_k)\in (\C^\ast)^{k-1}$ and $S\in \max \C[t_2^{\pm 1},\cdots,t_k^{\pm 1}]$. Let $\bar\Xi$ be the set of orbits in $\Xi$ under this action and given $\xi\in \Xi$, let $\bar\xi$ denote the orbit of $\xi$ in $\bar\Xi$.

 Define a map $\bchi:\Pi\rightarrow \Xi$ as follows.
  For $\pi\in \Pi$, let $$\bchi(\pi)(M) =
  \left\lbrace\begin{array}{ll} (\pi(M)(K_1), \pi(M)|_{\lie h_\fn}\mod Q_{fin}), &\text{if}\ M\in \supp(\pi)\\
  (0,0), & \text{otherwise}\end{array}\right.$$
\label{map.X}

\rems (1) It follows from \corref{weyl} that if $X_{\psi}$ is a $\cal L^c(\lie g)$-irreducible constituent of $W_\pi$ for $\pi, \psi\in \Pi$, then $\bchi(\pi) = \bchi(\psi)$. 

\noindent (2) Since $L(W_\pi)$ is spanned by elements of the form $w\otimes f$ with $w\in W_\pi$ and $f\in \C[t_2^{\pm1},\cdots,t_k^{\pm1}]$,  it follows from \corref{weyl} and \propref{gr.weyl} that if $X_\psi^\bos$ is $\cal T_k(\lie g)$-irreducible constituent of $L^\bos(W_\pi)$, then $\bchi(\psi)=\bchi(\pi).$

\subsection{} As $\cal J_{int}^+$ is an abelian category consisting of finite-length objects of $\cal I_{fin}^{(m\boe_1)}$, for each $m\in \mathbb N$, the objects in $\cal J_{int}^m$ can be written as the direct sum of indecomposables.
We say two indecomposable objects $U_1$ and $U_2$ in $\cal J_{int}^+$ are linked and write $U_1\sim U_2$ if there exists a family of indecomposable objects $U_1=U$, $U_2, \cdots, U_r=V$ in $\cal J_{int}^+$ such that either $\hom_{\cal J_{int}^+}(U_i, U_{i+1})\neq 0$ or   $\hom_{\cal J_{int}^+}(U_{i+1}, U_i)\neq 0$ for all $i=1,\cdots r-1.$ It is easy to see that $\sim$ induces an 
equivalence  relation on $\cal J_{int}^+$ and the corresponding equivalence classes are called the blocks in $\cal J_{int}^+$. Each block is a full abelian subcategory and the category $\cal J_{int}^+$ is the direct sum of the blocks. Clearly, if $Ext^1_{\cal J_{int}^+}(X,X')\neq 0$ for two irreducible objects $X, X'$ in $\cal J_{int}^+$, then $X\sim X'.$ 


One of the ingredients that lead towards the block decomposition of $\cal J_{int}^+$ is the following vanishing result for $\Ext^1$. 


 \vspace{.25cm}
 
\begin{prop}\label{L(g)-link} Given $\pi, \psi\in \Pi$ such that $\overline{\bchi(\pi)}\neq \overline{\bchi(\psi)},$
$$\Ext^1_{\cal J_{int}^+}(X_{\pi}^{\bog}, X_{\psi}^{\bos}) =0,\hspace{.85cm} \forall \ \ {\bog}, {\bos}\in \Z^k_0.$$
In particular, if  $V$ is linked to $X_\pi^{{\bog}}$ and $\overline{\bchi(\pi)}\neq \overline{\bchi(\psi)}$ then 
$$\Ext^1_{\cal J_{int}^+}(V, X_{\psi}^{{\bos}}) =0,\hspace{.85cm} \forall \ \ {\bos}\in \Z^k_0.$$  
\end{prop}
The arguments in the proof are similar to those in \cite[Proposition 6.5]{CG1}. We give the details here for the sake of completeness.
\proof  For a contradiction assume that there exists an integrable indecomposable $\cal T_k(\lie g)$-module $V$ with finite-dimensional weight spaces such that
    \begin{equation} \label{seq}
      0 \rightarrow X_{\psi}^{\bos} \overset{\i}{\rightarrow} V \overset{\p}{\rightarrow} X_{\pi}^{\bog}\rightarrow 0
      \end{equation}
    is a non-split short exact sequence. Applying the exact contravariant functor $\tilde{\omega}$ (if necessary) we may assume that $\wt(\psi)-\wt(\pi)\notin Q^+$.

Let $\wt(\pi)=\Lambda_\pi$ and $\wt(\psi) = \Lambda_\psi$. Since $V$ is a finite-length object of $\cal J_{int}^+$, by \corref{V-aff+}, $V_{aff}^+$ is a non-zero subspace of $V$. Let $v\in V_{aff}^+$ be a non-zero vector, then $\p(v)\in X_{\pi}^\bog$ is such that $X.\p(v) =0$ for all $X\in \lie n_{aff}^+\otimes \C[t_2^{\pm1},\cdots,t_k^{\pm1}]$. As $\Lambda_\psi\not> \Lambda_\pi$, either  $\p(v)\in (X_{\pi}^{\bog})_{aff}^+$ or $\p(v)=0$. If $\p(v)\in (X_{\pi}^{\bog})_{aff}^+$, then the weight of $\p(v)$ and hence of $v$ with respect to $\lie h_{tor}$, is $\Lambda_{\pi}+\delta_{\bop}$, where $\bop\in \Z^k_0$ is such that $\bop-\bog\in G_\pi$. If $\p(v)=0$, then $v\in \ker(\p) = \Im (\i)$. Hence, $v$ is the image under $\i$ of an element in $(X_\psi^\bos)_{aff}^+$ and the weight of $v$ with respect to $\lie h_{tor}$ is $\Lambda_{\psi}+\delta_{\bor}$, where $\bor\in \Z^k_0$ is such that $\bor-\bos\in G_\psi$. 

Since \eqref{seq} is non-split and $\Lambda_\psi\not> \Lambda_\pi$, there 
exists a non-zero vector $v_0\in V_{aff}^+$ such that 
$\p(v_0) = v_{\pi}\otimes t^{\bog}$. Let $\cal V = \cal U(\cal T_k(\lie g)).v_0$ and $\cal V_0 = \cal U(\cal T_{k-1}(\lie h_{aff})).v_0$. As $\Lambda_\psi\not> \Lambda_\pi$, for all $v\in \cal V_0$, $X.\p(v) =0$ for $X\in \lie n_{aff}^+\otimes \C[t_2^{\pm1},\cdots,t_k^{\pm1}].$ Hence it follows from the above discussion that $\cal V_0\subset \cal V_{aff}^+ \subseteq V_{aff}^+$. Furthermore, if $v_1,v_2\in \cal V_0$ is such that $\p(v_1)=\p(v_2)$, then $v_1-v_2\in (X_{\psi}^\bos)_{aff}^+$. Thus we get,
\begin{equation} \label{V.0}
\cal V_0 = \underset{\bom\in \Z^k_0}{\oplus} \cal V_{\Lambda_\pi +\delta_{\bom}} \cong_{\cal T_{k-1}(\lie h_{aff})} 
\underset{\bom\in \Z^k_0}{\oplus} \ ( (X_{\pi}^{\bog})_{\Lambda_\pi+\delta_{\bom}}+(X_{\psi}^{\bos})_{\Lambda_\pi+\delta_{\bom}}). \end{equation}
If $\dim \cal V_{\Lambda_\pi+\delta_{\bom}}\leq 1$ for all $\bom\in\Z^k_0$, then 
$\cal V$ is a highest $\ell$-weight module and, hence, by \propref{gr.weyl}, it 
is a quotient of $L^\bog(W_{\pi})$. This implies that every irreducible quotient of $\cal V$ is also an irreducible constituent of $L^\bog(W_\pi)$. But, by 
\remref{L(g)-link}(2), if $X_\psi^\bos$ is a $\cal T_k(\lie g)$-irreducible 
constituent of $L^\bog(W_\pi)$, then $\overline{\bchi(\pi)}=\overline{\bchi(\psi)}$. Hence, $\cal V$ cannot be a highest $\ell$-weight module. 

Thus,  there exists $\bor\in\Z^{k}_0$ such that 
$\dim (\cal V_0)_{\Lambda_\pi+\delta_{{\bor}}}>1.$ By \eqref{V.0}, this is possible 
only if $\Lambda_\pi=\Lambda_\psi=\Lambda$ and, since every irreducible 
$\cal T_k(\lie g)$-module is a highest $\ell$-weight module, 
$\dim \cal V_{\Lambda+\delta_{{\bor}}} =2$.  In particular, 
$\dim \cal V_{\Lambda+\delta_{\bom}} =2$, for all  $\bom\in \Z^{k}_0$ such that  
${\bom}- {\bor}\in G_{\pi}\cap G_{\psi}$. For each $\bom\in G_{\pi}\cap G_{\psi}$, fix an ordered basis $\{u_{1,\bom}, u_{2,\bom}\}$ of $\cal V_{\Lambda+\delta_{\bom}}$ such that $\p(u_1,\bom)\in X_\pi^\bog$ and $u_{2,\bom}\in \i(X_\psi^\bos)$, that is, 
\begin{equation} \label{p.i.def}\p(u_1,\bom) = v_\pi\otimes t^\bom, \hspace{.25cm}\text{and}\hspace{.25cm} u_{2,\bom} = \i(v_\psi\otimes t^\bom).\end{equation} 
For $i=1,2$, define algebra homomorphisms
  $\phi_i: \cal U(\cal T_{k-1}(\lie h_{aff}))\rightarrow \C$, by
  $$\phi_1(h\otimes t^{{\bon}}) = \sum_{M\in \supp(\pi)} \pi(M)(h)\ev_M(t^{{\bon}}), \hspace{.5cm} \forall\ \ h\otimes t^{{\bon}}\in \cal T_{k-1}(\lie h_{aff}),$$  
$$\phi_2(h\otimes t^{{\bon}}) = \sum_{M\in \supp(\psi)} \psi(M)(h)\ev_M(t^{{\bon}}), \hspace{.5cm} \forall\ \ h\otimes t^{{\bon}}\in \cal T_{k-1}(\lie h_{aff}).$$
Since $\i$ is injective, by \eqref{T(h).action} and \eqref{p.i.def}, for all $\bon\in G_{\psi}$, $h\otimes t^\bon.u_{2,\bom} =\phi_2(h\otimes t^\bon)u_{2,\bom+\bon}.$ Also, for any $c\in \C$, the vector $u_{1,\bom}^c= u_{1,\bom}+cu_{2,\bom}\in \cal V_0$ is such that $\p(u_{1,\bom}^c) = v_{\pi}\otimes t^\bom$. Therefore, for $\bon\in G_\pi\cap G_\psi$, $h\otimes t^\bon.u_{1,\bom}$ lies in the span of $\{u_{1,\bom+\bon}, u_{2,\bom+\bon}\}$ and is of the form  
$h\otimes t^\bon.u_{1,\bom} = \phi(h\otimes t^\bon)u_{1,\bom+\bon} + C u_{2,\bom+\bon}$. 

 Since $\overline{\bchi(\pi)}\neq \overline{\bchi(\psi)}$, $\pi$ and $\psi$ are distinct functions 
in $\Pi$. Hence by \eqref{T(h).action}, for some 
$\bom \in G_{\pi}\cap G_{\psi}$, there exists  
$H\in \cal U(\cal T_{k-1}(\lie h_{aff} ))_{\bon}$ such that 
$\phi_1(H)\neq \phi_2(H)$ and 
\begin{equation}\label{H.u.1.2} H.u_{1,{\bom}} = \phi_1(H)u_{1,{\bom}+{\bon}} + c u_{2,{\bom}+{\bon}}, \hspace{1cm} 
H.u_{2,\bom} = \phi_2(H)u_{2,\bom+\bon},\end{equation} where $c\in \C.$ Let $\varrho: \cal V_{\Lambda+\delta_{\bom+\bon}} \rightarrow  \cal V_{\Lambda+\bom}$ be the isomorphism of vector spaces defined by
  $\varrho(u_{i,\bom+\bon}) =u_{i,\bom}$ for $i=1,2$. Then
  the matrix of $\varrho \circ H: \cal V_{\Lambda+\delta_\bom} \rightarrow \cal V_{\Lambda+\delta_\bom}$ with respect to the ordered basis $\{u_{2,\bom},u_{1,\bom}\}$ is upper triangular and has  distinct diagonal entries namely $\phi_2(H)$ and
  $\phi_1(H)$. This implies that $\rho\circ H$ is diagonalizable. Since 
  $u_{2,\bom}$ is an eigenvector with  eigenvalue $\phi_2(H)$, for some $c\in \C$, $u_{1,\bom}^c$ must be an eigenvector of $\varrho\circ H$ corresponding to the eigenvalue $\phi_1(H)$. On the
  other hand as $u_{2,\bom} \in \cal U(\cal T_{k-1}(\lie h_{aff})).u^c_{1,\bom}$, there exists $H_0\in \cal U(\cal T_{k-1}(\lie h_{aff}))_{\bold{0}}$ such that \begin{equation} \label{u.2.1} u_{2,\bom} = H_0u^c_{1,\bom}.\end{equation} Since $\overline{\bchi(\pi)}\ne \overline{\bchi(\psi)}$,  one of the following holds.
\begin{enumerate}
\item[(i)] For all $\bob\in (\C^\ast)^{k-1}$, $\supp (\pi) \ne \bob.\supp(\psi)$. 
\item[(ii)] $\supp(\pi) =\bob.\supp(\psi)$ and there exists 
$M\in \supp(\pi)\cap \bob.\supp(\psi)$ such that $\psi(M)\notin \pi(M)+Q_{aff}$. 
\end{enumerate}

Suppose (i) holds. Then for each $\bob\in (\C^\ast)^{k-1}$, there exists $\lie m \in \supp(\pi)$ such that $\lie m\notin \bob.\supp(\pi)$. 
Let 
$$ I_{\lie m} = (\underset{N\in \bob.\supp(\psi)}{\bigcap} N)\cap (\underset{N\in \supp(\pi)-\{\lie m\}}{\bigcap} N) .$$ Then choosing $f\in I_{\lie m}$,  we see that $h\otimes f.u_{2,\bom} = 0$ but $h\otimes f.H_0u_{1,\bom} \ne 0$ which is a contradiction to \eqref{u.2.1}. 

Now suppose (ii) holds and 
$$I_M = (\underset{N\in \bob.\supp(\psi)-\{M\}}{\bigcap} N) .$$
 Then, by \eqref{T(h).action} and \eqref{H.u.1.2},  for $f  \in I_ M$ of the form $f = \sum_{i=1}^l a_it^{\bon_i} $, 
\begin{equation} \label{heq1} h\otimes f.u_{2,\bom} = \psi(M)(h) (\sum_{i=1}^l a_i\ev_M(t^{\bon_i}) u_{2, \bon_i+\bom}).\end{equation}  On the other hand, as $[h\otimes f, H_0] =0$, we get
  \begin{equation} h\otimes f.H_0.u^c_{1,\bom} = (\pi(M)(h) H_0.(\sum_{i=1}^l a_i \ev_M(t^\bon_i) u_{1, \bon_i+\bom}) + X, \label{heq2}\end{equation} where $X$ lies in the linear span of $\{u_{2, \bom +\bon_i} : i=1,\cdots, l\}$. Since $f\ne 0$(i,e, there exists $1\le i\le l$ such that $a_i\ne 0$), $\pi(M)$ is a dominant integral weight, for each $\bor \in G_\pi\cap G_\psi$, $\{u_{1,\bor}, u_{2,\bor}\}$ is linearly independent and $\ev_M{t^\bom}\neq 0$ for $M\in \max\C[t_2^{\pm1},\cdots,t_k^{\pm1}]$, it follows from \eqref{heq1}, \eqref{heq2} and \eqref{u.2.1},  that there exists $1\le i\le l$ for which $H_0.u_{1,\bon_i+\bom} =0$. This imples that $H_0.v =0$ for all $v\in \cal V_0$. In particular, $u_{2,\bom}=0$ which is a contradiction. Hence the proposition.  

As every object in $\cal J_{int}^+$ has finite-length, applying induction on the length of a module in $\cal J_{int}^+,$  it is easy to see that, if $V$ is an object of $\cal J_{int}^+$ which is linked to $X_\pi^\bog$, then $\Ext_{\cal J_{int}^+}(V,X_\psi^\bos) =0$ whenever $\overline{\bchi(\psi)}\ne \overline{\bchi(\pi)}$.\endproof

\subsection{}  We now define the finitely supported functions of type $\bi$ and type $\bi\bi$.  
\label{diagram}
\begin{defn} 
 Given a finitely supported function $\pi\in \Pi$ we shall say  \begin{itemize}
\item[(i).] $\pi$ is of type $\bi$ if for all $M\in \supp\pi$,  $\pi(M)(K_1)<(\theta|\theta)$.
\item[(ii).] $\pi$ is of type $\bi\bi$ if there exists $M\in \supp \pi$ such that $\pi(M)(K_1)\geq (\theta|\theta)$.
\end{itemize}
\end{defn}
By abuse of language we shall refer to an irreducible $\cal T_k(\lie g)$-module $X_\pi^\bor$ as type $\bi$(respectively type $\bi\bi$) module whenever $\pi\in \Pi$ is of type $\bi$(respectively type $\bi\bi$).

\vspace{.15cm}
\label{xi.type}

\begin{prop} Let $\pi\in \Pi$ be a finitely supported function of type $\bi$. 
Then for $\bog\in \Z^k_0$, $L^\bog(W_\pi)$ is an irreducible 
$\cal T_k(\lie g)$-module. \end{prop}
\proof  Let $\pi_{\Lambda,M} \in \Pi$ be the finitely supported function such 
that  $$\supp(\pi_{\Lambda,M})=\{M\}\hspace{.5cm} \text{and} 
\hspace{.5cm}\pi_{\Lambda,M}(M) =\Lambda.$$
By \propref{weyl}(iv), we have $$W_\pi \cong_{\cal L^c(\lie g)} \underset{M\in \supp(\pi)}{\otimes} W_{\pi_{\pi(M),M}}.$$ We show that when $\pi\in \Pi$ is of type $\bi$, then for each $M\in \supp(\pi)$, $W_{\pi_{\pi(M),M}}$ is irreducible as a $\cal L^c(\lie g)$-module and hence is isomorphic to its irreducible quotient $X_{\pi_{\pi(M),M}}$.  As a consequence, $W_\pi\cong_{\cal L^c(\lie g)} X_\pi$. By \propref{weyl}(ii) and \propref{gr.weyl}, it thus follows that in this case, $L^\bog(W_\pi) \cong_{\cal T_k(\lie g)} X_\pi^\bog$ for all $\bog\in \Z^k_0$ and hence the result.

Since the $\cal L^c(\lie g)$-module $W_\pi$ has finite-dimensional weight 
spaces,  and $\pi(M)(\alpha_{n+1}^\vee) >0$ for all $M\in \supp(\pi)$, by 
\propref{affine}(ii), $W_{\pi_{\pi(M),M}}$ is completely irreducible when considered
as a module for $\lie g_{aff}$. 

For each $\bom\in \Z^k_0$, consider the map $\phi_{\bom}: \lie g_{aff}t^{\bom}\otimes W_{\pi_{\pi(M),M}} \rightarrow W_{\pi_{\pi(M),M}}$ defined by $$\phi_{\bom}(x\otimes t^{\bom},w)=(x\otimes t^{\bom})w, \hspace{.35cm} \text{for}\ x\in \lie g_{aff},\ w\in W_{\pi_{\pi(M),M}}.$$ Clearly, $\phi_{\bom}$ is a $\lie g_{aff}$-module map, where  the the action of  $\lie g_{aff}$ on the first factor is given by the adjoint representation. Now 
replacing the maps $\{\phi_r\}_{r\in \Z}$ in \cite[Theorem 4]{CL} by the maps $\{\phi_{\bom}\}_{{\bom}\in \Z^{k}_0}$  and using the result \propref{affine}(iv) in place of \cite[Theorem 3]{CL}, the same proof as \cite[Theorem 4]{CL} shows that $W_{\pi_{\pi(M),M}}$ is an irreducible $\cal L^c(\lie g)$-module whenever $\pi(M)(\alpha_{n+1}^\vee) < \theta(\theta^\vee).$ This completes the proof of the proposition.
 \endproof


\subsection{}\label{xi.1} The following is a result on the vanishing of $\Ext^1$ betweem two type $\bi$ irreducible $\cal T_k(\lie g)$-modules.
\begin{prop} Suppose $\lie g_{fin}$ is not of type $B_n$, $C_n$, $F_4$ or $G_2$. Let $\pi,\psi \in \Pi$ be two finitely supported functions of type $\bi$ such that $\overline{\bchi(\pi)} = \overline{\bchi(\psi)}$. Then $\Ext^1_{\cal J_{int}^+}(X_\pi^{\bog},X_\psi^{\bos})=0$ if $X_{\pi}^\bog$ is not isomorphic to $X_{\psi}^{\bos}$ as 
$\cal T_k(\lie g)$-modules. \end{prop}
\proof Firstly note that, $\theta(\theta^\vee) =4$ when $\lie g_{fin}$ of type $C_n$, and $\theta(\theta^\vee) =2$ otherwise. Therefore, if $\lie g_{fin}$ is not of type $C_n$ and  $\pi\in \Pi$ is of type $\bi$, then for all $M\in \supp(\pi),$ $\pi(M)(K_1) = 1.$
Now observe that if $\lie g_{fin}$ is not of type $B_n$, $C_n$, $F_4$ or $G_2$, then by 
\eqref{J.0}, \eqref{Lambda.i} and Table \ref{table:affine}, given 
$\gamma\in \Gamma$, there exists a unique $\Lambda\in P_{aff}^+$ such that 
$\Lambda(K_1) =1$ and $\Lambda|_{\lie h_{fin}} \equiv \gamma \mod 
Q_{fin}$. Hence, if $\supp(\pi)= \bob.supp(\psi)$, for $\bob\in (\C^\ast)^{k-1}$ then 
for each $M\in \supp(\pi)\cap \bob.\supp(\psi)$, $\pi(M) = \psi(M_\bob)$, where $M_\bob\in \supp(\psi)$ is such that 
$\bob.M_{\bob} = M$. By \corref{isomorphism} it thus follows that, $X_\pi$ 
is isomorphic to $X_\psi$ as a $\cal L^c(\lie g)$-module. Thus to prove the result
it is sufficient to show that $\Ext^1_{\cal J_{int}^+}(X_\pi^{\bog},X_\pi^{\bos})=0$ 
whenever $\bog-\bos\notin G_\pi$.

For a contradiction assume that, given $\bog,\bos\in \Z^k_0$ with $\bog-\bos\notin G_\pi$, there exists an indecomposable $\cal T_k(\lie g)$-module $V$ such that 
\begin{equation}\label{seq.3}
0\rightarrow X_{\pi}^{\bos}\overset{\i}{\rightarrow} V\overset{\p}{\rightarrow} X_{\pi}^{\bog}\rightarrow 0,
 \end{equation} is a non-split short exact sequence in $\cal J_{fin}^+.$ Since for all $\bog\in \Z^{k}_0$,  
$X_\pi^\bog \cong_{\cal T_k(\lie g)} X_{\pi}^{\bold{0}}\otimes \C_{\delta_{\bog}}$, tensoring the sequence \eqref{seq.3} by $\C_{\delta_{-\bor}}$ for some $\bor\in \Z_0^k$ if necessary, we may assume that $\bog\in G_\pi$. 

As $P_{aff}(X_{\pi}^{\bos})\subseteq \wt(\pi)-Q_{aff}^+$, there exists a non-zero weight vector $v\in V$ such that $\p(v) = v_\pi\otimes t^{\bog}$ and 
$$\begin{array}{ll}\lie n_{aff}^+\otimes \C[t_2^{\pm 1},\cdots,t_k^{\pm 1}].v=0, \hspace{.5cm} & h.v = \sum_{M\in \supp(\pi)}\pi(M)(h).v, \hspace{.25cm} 
\forall \ h\in \lie h_{aff}, \\ & \\
 h\otimes f.v=0, \hspace{.25cm} \forall\  h\in \lie h_{aff}, \ f \in \underset{M\in supp(\pi)}{\cap} M, \hspace{.5cm} & c.v=0 ,\ \forall c\in \cal Z_1'.\end{array}$$ 
Let $\cal V = \cal U(\cal T_k(\lie g)).v.$ Then, following the same arguments 
as in \propref{xi.type}, we see that $\cal V_{aff}^+$ is a non-zero subspace of 
$V_{aff}^+$. Furthermore, if $\cal V_0 = \cal U(\cal T_{k-1}(\lie h_{aff})).v$, 
then,  $$(\cal V_0)_{\Lambda_\pi +\delta_{\bom}} \cong_{\cal T_{k-1}(\lie h_{aff})} (X_\pi^\bog)_{\Lambda_\pi+\delta_{\bom}} \oplus (X_\pi^\bos)_{\Lambda_\pi+\delta_{\bom}}, \hspace{.45cm} \forall \ \bom\in \Z^k_0. $$
By \propref{six.five}, $\Lambda_{\pi}+\delta_{\bom}\in P(X_\pi^\bog)$ if and only if $\bog-\bom\in G_\pi$. Since, by choice $\bog-\bos\notin G_\pi$, given $\bom\in \Z^k_0$ such that $(\cal V_0)_{\Lambda_\pi+\delta_\bom}\ne 0$, 
$$(\cal V_0)_{\Lambda_\pi +\delta_{\bom}} \cong_{\cal T_{k-1}(\lie h_{aff})} (X_\pi^\bog)_{\Lambda_\pi+\delta_{\bom}} \hspace{.25cm} \text{or} \hspace{.25cm} (\cal V_0)_{\Lambda_\pi +\delta_{\bom}} \cong_{\cal T_{k-1}(\lie h_{aff})} (X_\pi^\bos)_{\Lambda_\pi+\delta_{\bom}}. $$
This implies that $\dim (V_0)_{\Lambda_\pi+\delta_\bom}\leq 1$ for each $\bom\in \Z^k_0$, that is, $\cal V$ is a highest $\ell$-weight module with unique irreducible quotient $X_\pi^\bog$. Hence by \propref{gr.weyl}, $\cal V$ is a quotient of $\L^\bog(W_\pi)$. But by \propref{xi.type}, $L^\bog(W_\pi)$ is an irreducible $\cal T_{k}(\lie g)$-module whenever $\pi$ is of type $\bi$. So, $\cal V\cong_{\cal T_k(\lie g)}L^\bog(W_\pi) \cong_{\cal T_k(\lie g)} X_\pi^\bog$ and by \propref{isomorphism}, 
$X_\pi^{\bos}\cap \cal V= 0$. This proves the proposition.
\endproof

\vspace{.15cm}

\rem There are two reasons why in the above proposition we exclude the cases when $\lie g_{fin}$ is of 
type  $B_n$, $C_n$, $F_4$ and $G_2$. Firstly, for $\lie g_{fin}$ of type $B_n$, $\mu_1= \Lambda_{n+1}$ and $\mu_2= \Lambda_{n+1} +\omega_1$ are both dominant integral weights for which $\mu_1(K_1) = \mu_2(K_1) =1$ and $(\mu_2 - \mu_1)|_{\lie h_{fin}} \in Q_{fin}$. Hence in this case,  the arguments used in proposition do not work. Moreover, since $\mu_i(K_1)<\theta(\theta^\vee)$, the methods used in \secref{hom}  to show the linkage between two irreducible $\cal T_k(\lie g)$-modules of type $\bi\bi$ cannot be used.  Similarly, for $\lie g_{fin}$ of type $G_2$(respectively, type $F_4$), $\mu_1= \Lambda_{n+1}$ and $\mu_2= \Lambda_{n+1} +\omega_2$ (respectively, $\mu_2 = \Lambda_{n+1} +\omega_4$) are both dominant integral weights for which $\mu_1(K_1) = \mu_2(K_1) =1$ and $(\mu_2 - \mu_1)|_{\lie h_{fin}} \in Q_{fin}$ and in the case when $\lie g_{fin}$ is of type $C_n$, for 
$l=1,2,3$ there exists distinct dominant integral weights $\mu, \mu'$ such that $\mu(K_1) = \mu'(K_1) = l$ and $(\mu - \mu')|_{\lie h_{fin}} \in Q_{fin}$, so the arguments of \propref{xi.1} do not work. Also, since, in this case, $\theta(\theta^\vee) =4$, the methods of \secref{hom}) to show the linkage between irreducible $\cal T_k(\lie g)$-modules of type $\bi\bi$ do not work.

\section{The Main Results}\label{Proof.Prop}

\subsection{} We begin with the following proposition.
\label{hom}
\begin{lem} Let $\mu_1, \mu_2\in P_{aff}^+$. Supppose there exists a non-zero $\lie g_{aff}$-module homomorphism $\phi: \lie g_{aff}\otimes X(\mu_1)\rightarrow X(\mu_2)$. Then, for each $M\in \max \C[t_2^{\pm1},\cdots,t_k^{\pm1}]$, the following formulas define an action of $\cal L^c(\lie g)$-module on $X(\mu_1)\oplus X(\mu_2)$:
\begin{equation}\label{eq.af.st} x\otimes t^{{\bom}}(v,w) = (\ev_M(t^{{\bom}})x.v, \ev_M(t^{{\bom}})x.w+ \ev_M(\underset{i=2}{\overset{k}{\sum}} \frac{\partial}{\partial t_i} (t^{{\bom}}) )\phi(x\otimes v)),
\end{equation} where $x\in \lie g_{aff}$, $\bom\in \Z^k_0$, $v\in X(\mu_1)$, $w\in X(\mu_2)$. Denoting this module by $X(\mu_1,\mu_2,M)$, it follows that
\begin{equation}\label{seq.12} 0\rightarrow X_{\pi_{\mu_2,M}} \rightarrow X(\mu_1,\mu_2,M) \rightarrow  X_{\pi_{\mu_1,M}} \rightarrow 0 \end{equation} 
is a non-split short exact sequence of $\cal L^c(\lie g)$-modules. In particular if $\mu_1 > \mu_2$, there exists a canonial $\cal L^c(\lie g)$-module surjective homomorphism $W_{\pi_{\mu_1,M}} \rightarrow X(\mu_1,\mu_2,M).$
\end{lem}
 \proof It is straightforward to check that the formula \eqref{eq.af.st} gives a $\cal L^c(\lie g)$ module structure on $X(\mu_1)\oplus X(\mu_2)$.
As $\cal L^c(\lie g).X_{\pi_{\mu_2, M}} \subseteq X_{\pi_{\mu_2,M}}$, $ X_{\pi_{\mu_2,M}}$ is a $\cal L^c(\lie g)$-submodule of $X(\mu_1,\mu_2,M)$. Furthermore, since $\phi$ is a non-zero homomorphism we see that $X(\mu_1,\mu_2,M)$ is an indecomposable $\cal L^c(\lie g)$-module and the sequence \eqref{seq.12} is non-split.  For the second part of the proposition, the proof is the same as in \cite[Proposition 3.4]{CM}. \endproof 

\vspace{.1cm}

\begin{prop} Let $\pi_1, \pi_2 \in \Pi$ be two finitely supported functions of type $\bi\bi$ such that $\overline{\bchi(\pi_1)} = \overline{\bchi(\pi_2)}$.  
Suppose 
$\lie g_{fin}$ is of type $A_n, D_n, E_6, E_7, E_8$.
Then, there exists $\tilde{\pi}_1\in \Pi$  and a sequence $\psi_0,\psi_1,\cdots,\psi_r$ of finitely supported functions of type $\bi\bi$ with $\overline{\bchi(\psi_j)} = \overline{\bchi(\pi_1)}$
for $j=0,1,2,\cdots,r$ such that up to tensoring by one-dimensional modules $X_{\pi_1}$ is isomorphic to $X_{\tilde{\pi}_1}$, $\psi_0 = \tilde{\pi}_1$, $\psi_r=\pi_2$, and either $\Ext^1_{\cal L^c(\lie g)}(X_{\psi_j},X_{\psi_{j+1}})\ne 0$ or  $\Ext^1_{\cal L^c(\lie g)}(X_{\psi_j},X_{\psi_{j-1}})\ne 0$ for each $j\in \{0,1,\cdots,r\}$. 
 \end{prop}
\proof First consider the case when 
 $(\pi_1(M)|\Lambda_{n+1}) = (\pi_2(M)|\Lambda_{n+1})$, for all 
$M\in supp(\pi_1)$.  Since $\pi_1$ is a function of type $\bi\bi$, there exists $M\in \supp(\pi)$ such that $\pi(M)(K_1)\geq \theta(\theta^\vee)$. Let $\pi_1^1, \pi_1^2\in \Pi$ be such that $\supp(\pi_1^1) = M$, $\pi_1^1(M) = \pi_1(M)$, $\supp(\pi_1^2) = \supp(\pi_1)-\{M\}$, $\pi_1^2(N) = \pi_1(N)$ for all $N\in \supp(\pi_1)$. Then $\pi_1=\pi_1^1+\pi_1^2$.

By \propref{algo}, there exists $\Lambda\in P_{aff}^+$ with $\Lambda(K_1) = \pi_1(M)(K_1)$ such that $\Lambda = \pi_1(M)+\alpha$ for some $\alpha\in Q_{fin}$ and either 
$\Hom (\lie g_{aff}\otimes X(\pi_1(M)),X(\Lambda))\ne 0$ or $\Hom (\lie g_{aff}\otimes X(\Lambda),X(\pi_1(M)))\ne 0$. This implies that there exists a non-split sequence of the form
$$0\rightarrow X_{\pi_{\pi_1(M),M}} \overset{\i}{\rightarrow} V\overset{\p}{\rightarrow} X_{\pi_{\Lambda,M}}\rightarrow 0,$$ or $$0\rightarrow X_{\pi_{\Lambda,M}} \overset{\i}{\rightarrow} V\overset{\p}{\rightarrow} X_{\pi_{\pi_1(M),M}}\rightarrow 0.$$ Tensoring the non-split sequence by $X_{\pi_1^2}$ and setting $\psi_1 = \pi_{\Lambda,M}+\pi_1^2$, we see that either $\Ext_{\cal L ^c(\lie g)}^1(X_{\psi_1},X_{\pi_1})\ne 0$ or $\Ext_{\cal L ^c(\lie g)}^1(X_{\pi_1},X_{\psi_1})\ne 0$. Applying \propref{algo},  and repeating the method we get the desired sequence. 

Suppose there exists $M\in \supp(\pi_1)$ such that $(\pi_1(M)|\Lambda_{n+1}) \ne (\pi_2(M)|\Lambda_{n+1})$. Then defining $\tilde{\pi}_1\in \Pi$ by 
\begin{eqnarray*}
\tilde{\pi}_1(M) = \left\{\begin{array}{ll} \pi_1(M),  \text{if}\, M\in \supp(\pi_1)\, \text{ is such that}  \, ((\pi_2(M)|\Lambda_{n+1}) = (\pi_1(M)|\Lambda_{n+1})), \\ \pi_1(M) + ((\pi_2(M)|\Lambda_{n+1}) - (\pi_1(M)|\Lambda_{n+1})),\,  \text{if}\, M\in \supp(\pi_1)\, \text{ is such that}  \\ \hspace{7cm} (\pi_2(M)|\Lambda_{n+1}) \ne (\pi_1(M)|\Lambda_{n+1}), 
\end{array}\right.
\end{eqnarray*} we see that up to tensoring by one-dimensional modules, $X_{\tilde{\pi}_1}$ is 
isomorphic to $X_{\pi_1}$ by \corref{isomorphism}. Now, by the first part of the proof, there exists a sequence of functions $\psi_1,\cdots, \psi_r$ satisfying the desired conditions. This completes the proof of the proposition. 
\endproof


\vspace{.15cm}

\subsection{} For $l\in \N$, let $\bs^l=\{\Lambda\in P_{aff}^+ : \Lambda(K_1) = l, \ \text{and}\ \Lambda|_{\lie h_{fin}}(\theta^\vee) \leq l\}$.

\label{L.s.r.link.1}
\begin{prop} Let $\pi\in \Pi$ be a finitely supported function of type $\bi\bi$ such that $G_\pi$ is a proper subgroup of $\Z^k_0$. 
\begin{enumerate} 
\item[i.]  If for some $M\in \supp(\pi)$ with $\pi(M)(K_1)\geq \theta(\theta^\vee)$,  $\pi(M)-\alpha \in \bs^{\pi(M)(K_1)}$ for some $\alpha\in Q_{fin}^+$,   then
 there exists $\psi\in \Pi$ such that $L(X_\psi)$ is simple and $L(X_\psi) \sim L^\bos(X_\pi)$ for some $\bos\in \Z^k_0$.
\item[ii.] If for some $M\in \supp(\pi)$ with $\pi(M)(K_1)\geq \theta(\theta^\vee)$, $\pi(M)+\alpha \in \bs^{\pi(M)(K_1)}$ for some $\alpha\in Q_{fin}^+$, then $L^\bos(X_\pi)\sim L^\bor(X_\pi)$ for all $\bos, \bor\in \Z_0^k$.
\end{enumerate}
\end{prop}
\proof Since $G_\pi$ is a proper subgroup of $\Z^{k-1}$, by  Equation \eqref{T(h).action} in \secref{six.five}, \begin{equation} \label{eq.0} \underset{N\in\supp(\pi) }\sum {\pi}(N)(\alpha_i^\vee)\ev_{N}(t^{\bom})=0, \hspace{.5cm} \forall \, i=1,\cdots,n+1 \hspace{.5cm} 
\forall\ \bom\in \Z^{k}_0-G_{\pi}. \end{equation} 
As $\pi\in \Pi$ is a function of type $\bi\bi$, by \propref{algo}, there exists $M\in supp(\pi)$ such that $\Lambda = \pi(M)+\beta\in \bs^{\pi(M)(K_1)}$ for some non-zero $\beta\in Q_{fin}$ and there exists a non-split sequence of the form
\begin{equation}\label{G.pi.seq} 0\rightarrow X_{\pi_{\pi(M),M}} \overset{\i}{\rightarrow} V\overset{\p}{\rightarrow} X_{\pi_{\Lambda,M}}\rightarrow 0,\hspace{.5cm} \text{or}\hspace{.5cm} 0\rightarrow X_{\pi_{\Lambda,M}} \overset{\i}{\rightarrow} V\overset{\p}{\rightarrow} X_{\pi_{\pi(M),M}}\rightarrow 0.\end{equation} 
Further, as $\beta\in Q_{fin}$ is non-zero, 
$\beta(\alpha_i^\vee)\ne 0$ for some $1\le i\le n$, and hence it will follow from \eqref{eq.0} that,
\begin{equation}\sum_{N\in supp(\pi)} {\pi}(N)(\alpha_i^\vee)\ev_{N}(t^{\bor}) + \beta(\alpha_i^\vee)\ev_{M}(t^\bor) \ne 0, \hspace{.5cm} \forall\ \,  \bor\in \Z_0^k.  \label{eq.1}\end{equation}
 
Let $\pi'\in \Pi$ be such that $\supp(\pi') = \supp(\pi)-\{M\}$ and $\pi'(N) = \pi(N)$ for all $N\in supp(\pi')$ and let $\psi\in \Pi$ be the function $\psi = \pi_{\Lambda, M} + \pi'$. Then by \eqref{eq.1}, $G_\psi = \Z^k_0$ and by \propref{gr.weyl}, $L(X_\psi)$ is an irreducible $\cal T_k(\lie g)$-module.

\item[i.] 
If $ \beta \in Q_{fin}^-$, then 
applying the functor $\tilde{\omega}$ if necessary, we get a non-split sequence 
\begin{equation} \label{seq.M.-}0\rightarrow X_{\pi_{\Lambda,M}} \overset{\i}{\rightarrow} V\overset{\p}{\rightarrow} X_{\pi_{\pi(M),M}}\rightarrow 0,\end{equation} with $\pi(M)|_{\lie h_{fin}} > \Lambda|_{\lie h_{fin}}$. 
Since tensoring by $X_{\pi'}$ is an exact functor and the functor $L$ preserves short exact sequences, it follows that 
\begin{equation} \label{seq.M.-}0\rightarrow L(X_{\psi}) \overset{\i}{\rightarrow} L(V\otimes X_{\pi'})\overset{\p}{\rightarrow} L(X_\pi)\rightarrow 0,\end{equation} is a short exact sequence. As $\wt(\pi)>\wt(\psi),$ and $X_\pi$ is an irreducible quotient of $V\otimes X_{\pi'}$, by \propref{gr.weyl}, $L(V\otimes X_{\pi'}) = \bigoplus\ L^\bos(V\otimes X_{\pi'})$ and for each $\bos\in \Z^k_0,$ 
$L^\bos(V\otimes X_{\pi'})$ is a quotient of $L^\bos(W_\pi)$. Thus, in this case there exists $\bos\in \Z^k_0$ such that $\Ext^1_{\cal J_{int}^+}(X_\pi^\bos, L(X_\psi))\ne 0$ and $\Ext^1_{\cal J_{int}^+}(X_\pi^\bor, L(X_\psi))= 0$ for all $\bor\in \Z_0^k$ such that $\bor-\bos\notin G_\pi.$

\item[ii.] If $ \beta \in Q_{fin}^+$, then 
applying the functor $\tilde{\omega}$ if necessary, we get a non-split sequence 
\begin{equation} \label{seq.M.+}0\rightarrow X_{\pi_{\pi(M), M}} \overset{\i}{\rightarrow} V\overset{\p}{\rightarrow} X_{\pi_{\Lambda,M}}\rightarrow 0,\end{equation} with $\pi(M)|_{\lie h_{fin}} < \Lambda|_{\lie h_{fin}}$. 
As above, tensoring by $X_{\pi'}$ and applying the functor $L$ to \eqref{seq.M.+} we obtain the short exact sequence  
\begin{equation} \label{seq.M.-}0\rightarrow L(X_{\pi}) = \oplus\ L^\bos(X_\pi) \overset{\i}{\rightarrow} L(V\otimes X_{\pi'})\overset{\p}{\rightarrow} L(X_\psi)\rightarrow 0.\end{equation} Since $\wt(\psi)>\wt(\pi),$ $X_\psi$ is an irreducible quotient of $V\otimes X_{\pi'}$, and $L(X_\psi)$ is an irreducible $\cal T_k(\lie g)$-module. Moreover, by \propref{gr.weyl}, $L(V\otimes X_{\pi'})$ is a highest $\ell$-weight quotient of $L(W_\pi)$ and for each $\bos\in \Z^k_0$, $L^\bos(X_\pi)$ is an irreducible submodule of $L(V\otimes X_{\pi'})$. This implies that, in this case,  the simple $\cal T_k(\lie g)$-modules $L^\bos(X_\pi)$ are linked to $L(X_\psi)$ for all $\bos\in \Z^k_0$ and completes the proof of the 
proposition. \endproof

\subsection{}  We say that a module $V\in \cal J^+_{int}$ has spectral character 
$\bar\xi\in \bar\Xi$ if $\overline{\bchi(\pi)} = \bar\xi$ for every
irreducible component $L^\bos(X_\pi)$ of $V$ . Given $\bar{\xi} \in \bar\Xi$, let $\cal J^+_{\bar{\xi}}$ be the subcategory consisting of all 
modules $V \in \cal J^+_{int}$ with spectral character $\bar\xi$.

 The following propsition proves that if $V, V'$ are irreducible objects of $\cal J^+_{\bar{\xi}}$ for some $\xi\in \Xi$, then $V$ is linked to $V'$. To prove the result we need the following elementary lemma. 

\label{eps}
\begin{lem} Let $\lambda_1,\cdots,\lambda_r$ be a set of dominant integral weights of $\lie g_{aff}$ and let $a_1,\cdots,a_r\in \C^\ast$ be a set of $r$ distinct non-zero complex numbers. Assume that there exists a positive integer $m$ such that for all $h\in \lie h_{aff}$, \begin{equation} \label{lem}\sum_{i=1}^r \lambda_i(h)a_i^s = 0, \hspace{.35cm} \text{whenever}\ s\not\equiv 0\mod m.\end{equation} Then $r\equiv 0\mod m$. Moreover there exists a permutation $\tau$ of $\{1,\cdots,r\}$ such that 
$$\begin{array}{llll}
\lambda_{\tau(1)}&=\lambda_{\tau(2)}&=\cdots&=\lambda_{\tau(m)},\\
\lambda_{\tau(m+1)}&=\lambda_{\tau(m+2)}&=\cdots&=\lambda_{\tau(2m)},\\
&\vdots & \vdots& \\
\lambda_{\tau(r-m+1)}&&=\cdots&=\lambda_{\tau(r)},
\end{array}$$
and complex numbers $a_{(1)},\cdots,a_{(p)}$ such that
$$\begin{array}{llll}
a_{\tau(1)}=\eps_m a_{(1)}, & a_{\tau(2)}=\eps_m^2 a_{(1)}, &\cdots, & a_{\tau(m)}=\eps_m^m a_{(1)},\\
a_{\tau(m+1)}=\eps_m a_{(2)}, & a_{\tau(m+2)}=\eps_m^2 a_{(2)}, &\cdots, & a_{\tau(2m)}=\eps_m^m a_{(2)
},\\ & \vdots &\vdots, & \\ a_{\tau(r-m+1)}=\eps_m a_{(p)}, & a_{\tau(r-m+2)}=\eps_m^2 a_{(p)},
&\cdots, & a_{\tau(r)}=\eps_m^m a_{(p)}, \end{array} $$
where $p=r/m$ and $\eps_m$ is a primitive $m^{th}$ root of unity.
\end{lem}

\begin{prop}  Let $\pi\in \Pi$ be a finitely supported function of type 
$\bi\bi$. If $\lie g_{fin}$ is of type $A_n, D_n, E_6, E_7$ or $ E_8$, 
then $X_\pi^\bos\sim X_\pi^\bor$ for all $\bor,\bos\in \Z^k_0$.
\end{prop}
\proof   Since $\pi$ is of type $\bi\bi$, there exists $M\in \supp(\pi)$ such 
that $\pi(M)(K_1)\geq \theta(\theta^\vee)$. If for some $M\in \supp(\pi)$ with $\pi(M)(K_1)\geq \theta(\theta^\vee)$, $\pi(M)+\alpha \in \bs^{\pi(M)(K_1)},$ for $\alpha\in Q_{fin}^+$, then by \propref{L.s.r.link.1},  $L^\bos(X_\pi)\sim L^\bor(X_\pi)$ for all $\bos, \bor\in \Z_0^k$.
It thus remains to prove the proposition in the case when $\pi(M)+\alpha\notin \bs^{\pi(M)(K_1)}$, for all $\alpha\in Q_{fin}^+$ and $M\in \supp(\pi)$ with $\pi(M)(K_1)\geq \theta(\theta^\vee)$.

If $G_\pi = \Z_0^{k}$, then $X_\pi^\bos = L(X_\pi) = X_\pi^\bor$ for all $\bos, \bor\in \Z^k_0$ and there is nothing to prove. 

Now assume that $G_\pi$ is a proper subgroup of $\Z_0^{k}$. By the structure theory of finitely generated abelian groups, there exists a basis $\boy_1, \boy_2,\cdots, \boy_{k-1}$ of $\Z^{k}_0$  such that $$G_\pi = m_1\boy_1\Z\oplus m_2\boy_2\Z\oplus \cdots \oplus m_{k-1}\boy_{k-1}\Z,$$ with  $m_i|m_{i+1}$ for $1\le i\le k-2$ and $m_i>1$ for at least one $i$. Then by  Equation \eqref{T(h).action} in \secref{six.five}, for all $h\in \lie h_{aff}$ and $1\le i\le k-1$,
\begin{equation}\label{relation.m_i}
\sum_{M\in \supp(\pi)} \pi(M)(h)\ev_M(\boy_i^s) =0 , \hspace{.35cm} \text{whenever}\ s\not\equiv 0\mod m_i. \end{equation}

We now prove the proposition by applying induction on the integer $k-1$. 

If $k-1 = 1$, then enumerating the maximal ideals in $\supp(\pi)$ appropriately, it will follow from  \lemref{eps} that  there exists maximal ideals $M_{(1)},\cdots,M_{(p)}$ of $\C[t_2^{\pm1}]$ such that 
\begin{equation}\begin{array}{lllll}
\pi(M_{1})&=\pi(M_2)&=\cdots&=\pi(M_m),& \cdots\\
\pi(M_{r-m+1})&&=\cdots&=\pi(M_r),&
\end{array}\label{M.ep.1}\end{equation}
and complex numbers $a_{(1)},\cdots,a_{(p)}$ such that
\begin{equation}\begin{array}{lllll}
M_1=\eps_m M_{(1)}, & M_2=\eps_m^2 M_{(1)}, &\cdots, & M_m=\eps_m^m M_{(1)},& \cdots,\\
M_{r-m+1}=\eps_m M_{(p)}, & M_{r-m+2}=\eps_m^2 M_{(p)},
&\cdots, & M_r=\eps_m^m M_{(p)},& \end{array} \label{M.ep.2}\end{equation}
where $p=r/m$ and $\eps_m$ is a primitive $m^{th}$ root of unity.

On the other hand, since $\pi$ is of type $\bi\bi$, and by the assumption, that for all $M\in \supp(\pi)$ with $\pi(M)(K_1)\geq \theta(\theta^\vee)$, $\pi(M)+\alpha \notin \bs^{\pi(M)(K_1)}$ for all $\alpha\in Q_{fin}^+$, there exists $N\in \supp(\pi)$ such that
$\pi(N)-\alpha \in \bs^{\pi(M)(K_1)}$ for some $\alpha\in Q_{fin}^+$.

Suppose $N=M_1$. Now let  
$j=1,\cdots,m$ set
$$\psi_j = \sum_{\underset{i\neq j}{i=1}}^r \pi_{\pi(M_i),M_i}+\pi_{\pi(M_j)-\alpha,M_j},$$ 
and let 
$$\psi = \sum_{i=m+1}^r \pi_{\pi(M_i),M_i}+\sum_{j=1}^m \pi_{\pi(M_j)-\alpha,M_j}.$$
Using \lemref{eps}, it is easy to see that $G_\pi = G_\psi = m_1\Z$  and by \propref{L.s.r.link.1}, $L(X_{\psi_j})$ is a simple $\cal T_2(\lie g)$-module such that $L(X_{\psi_j})$ is linked to $X_\psi^s$ for all $s\in \Z$. Furthermore, $\Ext^1_{\cal J_{int}^+}(X_\pi^r,L(X_{\psi_j})) \neq 0 $ for some $0\le r< m_1.$

Now using the conditions \eqref{M.ep.1} and \eqref{M.ep.2}, and the action of $h\otimes t_2$ on a highest weight vector of $L(X_{\psi_j})$, it can be easily verified using \eqref{T(h).action} that for $j=1,2,\cdots, m$, $\Z/m\Z$ acts by the irreducible character $\eps_{m_1}^{j-1}$ on $L(X_{\psi_j})$. Hence, following the arguments of \cite[Proposition 5.5(ii)]{CG1}, it follows that $L(X_{\psi_j})$ is linked to $X_\pi^{j-1}$. Since $L(X_{\psi_j})\sim X_{\pi}^{j-1}$ and $X_\psi^r \sim L(X_{\psi_j})$ for all $r, j\in \{0,1,\cdots,m_1\}$,  it follows from the transitivity of $\sim$ that $X_\pi^r\sim X_\pi^s$ for all $r,s\in Z$.

Now assume $k-1>1$. 
 Since $m_i|m_{i+1}$ for all $i$, and $m_i>1$ for some $i$, 
\begin{equation}\label{relation.m_i.2}
\sum_{M\in \supp(\pi)} \pi(M)(h)\ev_M((\boy_1\boy_2\cdots \boy_i)^s) =0 , \hspace{.35cm} \text{whenever}\ s\not\equiv 0\mod m_i. \end{equation} 
As $\pi(M)\in P_{aff}^+$ for all $M\in \supp(\pi)$ and \eqref{relation.m_i.2} holds for all $i\in \{1,\cdots,k-1\}$,  $\ev_M(\boy_1\boy_2\cdots \boy_i)$ cannot be equal for all $M\in \supp(\pi).$ Therefore grouping together the coeffients for all $M\in \supp(\pi)$ for which $\ev_M(\boy_1\boy_2\cdots \boy_i)$ are equal, we get an equation of the form \eqref{lem}.
Now using \lemref{eps} in \eqref{relation.m_i.2}, for $1\le i\le k-1$, and 
following the same arguments as above, it can be seen that $X^\bor_\pi \sim X^\bos_\pi$ for all $\bos, \bor\in \Z^k_0$.  \endproof

\subsection{} We finally state and prove the main result of the paper.

\begin{thm}  Let $\cal T_k(\lie g)$ be a toroidal Lie algebra with underlying finite-dimensional Lie algebra of one of the following types $A_n,D_n,E_6,E_7,E_8$. Then  every indecomposable $\cal T_k(\lie g)$-module $V$ in $\cal J_{int}^+$ is an object of $\cal J^+_{\bar\xi}$ for some $\xi\in \Xi$. Moreover, if  there exists $\pi\in \Pi$ of type $\bi$ such that $\overline{\bchi(\pi)} = \bar\xi$, then the irreducible components of $V$ are all isomorphic.
\end{thm}
\proof We prove the result by applying induction on the length of the indecomposable $\cal T_k(\lie g)$-module $V$.  Suppose $V$ is of length $1$. Then $V$ is irreducible and isomorphic to $X^\bos_\pi$ for some $\pi\in \Pi$. In this case $V\in \cal J^+_{\overline{\bchi(\pi)}}$ and we are done.

If $V$ is not irreducible then we have an extension 
\begin{equation} \label{pi.U} 0 \rightarrow X^\bos_\pi \overset{\i}{\rightarrow} V \overset{\p}{\rightarrow} U\rightarrow 0,\end{equation} for some $\pi\in \Pi$ and $\bos\in \Z^k_0$.

 As $ U\in \Ob \cal J_{int}^+$, it can be written as the direct sum of indecomposable $\cal T_k(\lie g)$-modules $U_j, j=1,\cdots,r$. Notice that the length of each $\cal T_k(\lie g)$-module $U_j$ is strictly less than the length of $V$. Therefore,  by inductive hypothesis there exists $\xi_j\in \Xi$ such that  every irreducible component $X_{\pi_{ij}}^{g_{ij}}$ of $U_j$ for $j=1,\cdots,r$ is such that $\overline{\bchi(\pi_{ij})}=\overline{\xi_j}$. If  
$\overline{\xi_j}\neq \overline{\bchi(\pi)}$ for some $1\leq j\leq r$, by \propref{L(g)-link}, $\Ext_{\cal J_{int}^1}(U_{j},X_\pi^\bog)=0$, which implies that there exists a direct summand of $V$ that is isomorphic to $U_{j}$. This contradicts our assumption that $V$ is indecomposable. 

If there exists $\pi\in \Pi$ of type $\bi$ such that $\overline{\bchi(\pi)} = \bar\xi$, then it follows from \propref{xi.1} that in the sequence \eqref{pi.U},
 every irreducible component of $U$ must be isomorphic to $X^\bos_\pi.$
\label{main}\endproof

\subsection{Remark} The  methods used in this paper cannot be extended to determine and characterize the blocks of $J_{fin}^+$ in the case when the underlying finite-finite Lie algebra of $\cal T_k(\lie g)$ is not simply-laced.  This is because the results proved here, are dependent on the fact that whenever $\lambda\in P_{fin}^+$ is such that $l\Lambda_{n+1}+\lambda \in P_{aff}^+$ for $l\geq \theta(\theta^\vee)$, there exists a sequence $\{\eta_{\lambda,r}\}_r$ of the form described in \propref{algo}, such that   $$ht^w(\eta_{\lambda,r}) \leq ht^w(\lambda)\le \lambda(\theta^\vee) \hspace{.25cm} \forall \ r\in \N,$$
and $\eta_{\lambda,r} = \omega_i$ for some $i\in J_0$ when  $\lambda\not\equiv 0\mod Q_{fin}$ and  $\eta_{\lambda,r} = \theta$ when  $\lambda\equiv 0\mod Q_{fin}$. 
However, for $\lie g_{fin}$ of type $B_n, C_n$, $F_4$ or  $G_2$, though  a sequence of the desired form can be obtained for all fundamental weights $\lambda$ with $ht^w(\lambda)>1$, one cannot obtain such a sequence for $\lambda\in P_{fin}^+$ with $ht^w(\lambda) =1$.

When $\lie g_{fin}$ is of type $B_n$, $C_n$ or $F_4$ and $\lambda$ is a fundamental weight with $ht^w(\lambda) =1$  and $\lambda\equiv 0\mod Q_{fin}$, we see that
$$ht^w(\lambda+w(\theta))> \max\{\lambda(\theta^\vee), ht^w(\lambda)\}, \hspace{.15cm} \text{and} \hspace{.15cm} ht^w(\theta+w(\theta)) \geq ht^w(\lambda) $$ for all $w\in W_{fin}$ with $\lambda+w(\theta), \theta+w(\theta) \in P_{fin}^+,$ and 
$$ht^w(\lambda+\theta-\beta) > \max\{\lambda(\theta^\vee), ht^w(\lambda)\}, \hspace{.15cm} \text{and} \hspace{.15cm} ht^w(2\theta-\beta) \geq ht^w(\lambda) $$
for all $\beta\in R_{fin}^+$ with $\lambda+\theta -\beta, 2\theta-\beta \in P_{fin}^+.$ 

For $\lie g_{fin}$ of type $C_n, n\geq 2$, $J_0 =\{1\}$, $\theta(\theta^\vee)$=4 and $ht^w(\lambda)=1$ for all fundamental weights $\lambda$. The case when $\lambda \equiv 0\mod Q_{fin}$ has been discussed above. Now, consider the case when $\lambda$ is a fundamental weight such that $\lambda \not\equiv 0\mod Q_{fin}$. Then, using the relations,
\begin{equation}\label{C.n} \begin{array}{c}
2\omega_i = 2\omega_{i-1}+ 2(\alpha_i+\alpha_{i+1}+\cdots+\alpha_{n-1})+\alpha_n, 
\hspace{.25cm} \text{for}\ 2\le i\le n, \hspace{.35cm} 2\omega_1 = \theta,\\
\omega_{i} = \omega_{1} +\omega_{i-1} + [\theta-(\sum_{i=1}^{i-1}\alpha_{j} +2(\alpha_{i}+\cdots+\alpha_{n-1})+\alpha_{n})], \hspace{.25cm} \text{for}\ 2\leq i\leq n,
 \end{array}
\end{equation} 
we see that any sequence $\{\eta_{j,\omega_i}\}_j$ for 
$2\le i\le n$, is such that $\eta_{j,\omega_i} = \omega_k+\omega_p$ for some $1\le k$ and $p=1$ or $2$ and $ht^w(\omega_k +\omega_p) > 1$. Hence, in this case also, one cannot obtain sequences of the  desired form.

Similarly, when $\lie g_{fin}$ is of type $G_2$, one cannot obtain a sequence of the desired form for $\lambda =\omega_2$. Using  \propref{tensor}(i) and (iii), we see that, \propref{tensor}(iii) cannot be applied here and for all $w\in W_{fin},$ $ht^w(\omega_2+w(\theta))> \max\{ht^w(\omega_2), \omega_2(\theta^\vee)\}$ and $ht^w(w(\omega_2)+\theta)> \max\{ht^w(\omega_2), \omega_2(\theta^\vee)\}$ whenever $\omega_2+w(\theta)\in P_{fin}^+$ or $w(\omega_2)+\theta\in P_{fin}^+$.

As a consequence, in each of the above cases,  one cannot obtain sequences of the  desired form for $\lambda = \sum_{i=1}^n r_i\omega_i \in P_{fin}^+$
with $r_i\neq 0$,  for $i\in I-J_0$ such that $ht^w(\omega_i) = 1.$


\begin{thebibliography}{}

\bibitem[A1]{A1}  D. ~Adamovi${\acute{\text{c}}}$, \emph{ Vertex operator algebras and irreducibility of certain modules for affine Lie algebras},~  Math. Res. Lett. \textbf{4} (1997), no. 6, 809-821. 

\bibitem[A2]{A2}  D.~ Adamovi$\acute{\text{c}}$, \emph{ An application of U(g)-bimodules to representation theory of affine Lie algebras}, ~ Algebr. Represent. Theory \textbf{7} (2004), no. 4, 457-469.
\bibitem[AL]{AL} J.~ Auger, M.~ Lau, \emph{Extension of modules for twisted current algebras}, ~arXiv:1704.03984 [math.RT].  
\bibitem[B]{Bou} N. Bourbaki, \emph{Elements of Mathematics, Lie Groups and Lie Algebras}, ~Chapters 4-6.

\bibitem[C]{C1} V.~Chari, \emph{ Integrable representations of affine Lie-algebras},~Invent. Math. \textbf{85} (1986), 317-335.

\bibitem[CFK]{CFK} V.~ Chari, G.~ Fourier, T.~ Khandai, \emph{A Categorical Approach to Weyl Modules}, Transformation Groups \textbf{15(3)}, (2010), 517-549.

\bibitem[CG]{CG1} V.~Chari and J.~Greenstein, {\em Graded Level Zero Integrable Representations of Affine Lie Algebras}, Transactions of AMS.
\textbf{360} (2008), no. 6, 2923-2940.

\bibitem[CL]{CL} V.~ Chari, T.~Le, \emph{Representations of Double Affine Lie Algebras},~ In a Tribute to C. S.~Seshadri (Chennai, 2002) Trends Math., Birkhauser, Basel, 2003, 199-219.

\bibitem[CM]{CM}  V. ~Chari, A. ~Moura, \emph{Spectral Characters of finite-dimensional representations of affine algebras}, ~J. Algebra \textbf{279} (2004), no.2, 820-839. 

\bibitem[CP1]{CPloop} V.~ Chari, A.~Pressley, \emph{New unitary representations of loop algebras}, ~Math. Ann. \textbf{275}, (1986), no.1 87-104.

\bibitem[CP2]{CP} V.~ Chari, A.~Pressley, \emph{A new family of Irreducible, Integrable Modules for Affine Lie algebras},~ Math. Ann. \textbf{277}, (1987), 543-562.

\bibitem[CP3]{CPweyl}V.~Chari and A.~Pressley, \emph{Weyl modules for classical andquantum affine  algebras}, Represent. Theory \textbf{5} (2001), 191--223 (electronic).




\bibitem[DGK]{DGK} V. V. Deodhar, O. Gabber, V. Kac \emph{Structure of some categories of representations of infinite-dimensional lie algebras} Advances in Mathematics \textbf{45}(1), (1982), 92-116.

\bibitem[FL] {FL} G. ~Fourier,  P. ~Littelmann, \emph{Weyl modules, Demazure
  modules,KR-modules, crystals, fusion products and limit constructions}, Adv.
  Math. \textbf{211}(2) (2007), 566–593.

\bibitem[H]{Humphreys} J.E.Humphreys, \emph{Introduction to Lie Algebras and Representation Theory}, Springer, 1968

 

\bibitem[Kac]{Kac} V. G.~Kac, \emph{Infinite dimensional Lie algebras} Prog. Math., Boston 44, Third Edition.

   
\bibitem[Kh]{Kh} T. ~Khandai, \emph{Irreducible Integrable Representations of Toroidal Lie Algebras}, arXiv:1607.04389

\bibitem[Ku]{icm.ku} S. ~Kumar, \emph{Tensor Product Decomposition}, Proceedings of the International Congress of Mathematicians, Hyderabad, India, 2010.
\bibitem[KL]{KL} D. ~Kus, P. ~Littelmann, \emph{Fusion products and toroidal algebras}, Pacific J. Math. \textbf{278} (2015) 427-445.
  
\bibitem[Ko]{Ko} R. ~Kodera, \emph{Extensions between finite-dimensional simple
  modules over a generalized current Lie algebra}, Transform. Groups
  \textbf{15} (2010), no. 2, 371-388.

\bibitem[MRY]{MRY} R.V.~ Moody, S. E.~ Rao, T.~ Yokomuma, \emph{Toroidal Lie algebra and vertex representations}, Geom.
Dedicata \textbf{35} (1990), 283-307.

\bibitem[MR]{MR} R.V.~Moody, S.E.~Rao \emph{Vertex Representations for N-Toroidal Lie algebras and a Generalization of the Virasoro Algebras}, Communications in Mathematical Physics \textbf{159} (1994), 239-264.

\bibitem[NS]{NS} E. ~ Neher, A. Savage, \emph{Extensions and block
  decompositions for finite-dimensional representations of equivariant map
  algebras}, Transform. Groups \textbf{20},  (2015), no. 1, 183-228.
  
\bibitem[S]{S} P. ~Senesi, \emph{The Block Decomposition of finite-dimensional representations of twisted loop algebras}, Pacific J. Math. \textbf{244} No. 2,(2010) 335-357.
  
\bibitem[R1]{R3} S. E. ~Rao, \emph{Complete reducibility of integrable modules for the affine Lie (super)algebras}, Journal of Algebra \textbf{264} (2003), no. 1, 269-278.

\bibitem[R2]{R2} ----, \emph{Classification of irreducible integrable modules for toroidal Lie algebras with finite dimensional weight spaces}, Journal of Algebra \textbf{277} (2004), 318-348.

  \bibitem[R3]{R05}----, \emph{Irreducible representations for toroidal Lie-algebras}, Journal of Pure and Applied Algebra \textbf{202}(1-3) (2005), 102-117.

    


\bibitem[RFS]{RSF}S. E.~Rao, V.~Futorny, S .~Sharma, \emph{Weyl Modules associated to Kac-Moody Lie Algebras}, Comm. Algebra \textbf{44} (2016), no. 12, 5045--5057.

\bibitem[VV]{VV} R.~ Venkatesh, S. ~Viswanath, Unique factorization of tensor products for Kac-Moody algebras. Adv. Math. \textbf{231} (2012), no. 6, 3162--3171

\end{thebibliography}
\end{document}